\documentclass{amsart}      

\usepackage{amsmath, amssymb, amsthm}
\usepackage{tikz, mathtools}
\usepackage{graphicx}
\usepackage{cite}

\def\Xint#1{\mathchoice
{\XXint\displaystyle\textstyle{#1}}%
{\XXint\textstyle\scriptstyle{#1}}%
{\XXint\scriptstyle\scriptscriptstyle{#1}}%
{\XXint\scriptscriptstyle\scriptscriptstyle{#1}}%
\!\int}
\def\XXint#1#2#3{{\setbox0=\hbox{$#1{#2#3}{\int}$ }
\vcenter{\hbox{$#2#3$ }}\kern-.6\wd0}}

\def\dashint{\Xint-}

\usepackage[width=5.8in,textheight=8.5in,letterpaper,centering]{geometry}

\newtheorem{theorem}{Theorem}
\newtheorem{lemma}{Lemma}
\newtheorem{proposition}{Proposition}
\newtheorem{remark}{Remark}

\begin{document}

\title{Chemical reaction-diffusion networks; convergence of the method of lines}

\author{Fatma Mohamed} 
\author{Casian Pantea} 
\author{Adrian Tudorascu}
\address{Department of Mathematics, West Virginia University, Morgantown, WV 26506}
\email{fmohamed@mix.wvu.edu, cpantea@math.wvu.edu, adriant@math.wvu.edu}
\thanks{AT was partially supported by  NSF Grant DMS-1600272 and by Grant \#246063 from the Simons Foundation. CP and FM  were partially supported by
NSF grant DMS-1517577.}

\maketitle

\begin{abstract}
We show that solutions of the chemical reaction-diffusion system associated to $A+B\rightleftharpoons C$ in one spatial dimension can be approximated in $L^2$ on any finite time interval by solutions of a space discretized ODE system which models the corresponding chemical reaction system replicated in the discretization subdomains where the concentrations are assumed spatially constant. Same-species reactions through the virtual boundaries of adjacent subdomains lead to diffusion in the vanishing limit. We show convergence of our numerical scheme by way of a consistency estimate, with features generalizable to reaction networks other than the one considered here, and to multiple space dimensions. In particular, the connection with the class of complex-balanced systems is briefly discussed here, and will be considered in future work.
\end{abstract}

\section{Introduction}
\label{intro}
Fueled in part by the advent of systems biology, the dynamical behavior of spatially homogeneous mass-action reaction systems has been the focus of much recent research.  A great number of results on the possibility of bistability or oscillation, local and global stability of equilibria, persistence of solutions etc. have been developed for ODE systems corresponding to well-mixed reaction networks.  This effort  started forty years ago \cite{Horn.1972aa, Horn.1972ab, Feinberg.1972aa}, and has seen a surge of interest in more recent years: \cite{Siegel.2000aa, Sontag.2001aa, Craciun.2005aa, Craciun.2006aa, Banaji.2007aa, Mincheva.2007aa, Angeli.2007aa, Conradi.2007aa, Craciun.2009aa, Anderson.2011sd, Craciun.2013ab, Banaji.2016aa}, to cite but a few examples. In particular, some of this work led to a proof of the Global Attractor Conjecture \cite{Craciun.gac}, a global asymptotic stability result for a large class of systems (called complex balanced networks).

On the other hand, much less is known about the corresponding reaction-diffusion setting, where the focus has largely been on the asymptotic behavior of solutions. One of the most studied examples is the reaction-diffusion system
$A + B{\leftrightharpoons} C$, whose solutions approach a spatially homogeneous distribution; this was shown by way of semigroup theory \cite{Rothe.1984aa}                                                                                                                                                                                                                                                        and entropy methods \cite{Desvillettes.2006aa}. Entropy considerations have also been used to successfully tackle other reaction-diffusion systems, including dimerization systems
$2A{\leftrightharpoons} B$ \cite{Desvillettes.2006aa}, weakly reversible monomolecular  reactions and other classes of linear systems \cite{Fellner.2015aa}, and classes of complex balanced systems with and without boundary equilibria \cite{Desvillettes.2016aa}. The latter work lays out a general method for complex balanced systems, but some of the technicalities depend on the specific network considered. This difficulty goes away under the assumption of equal diffusion coefficients, where general results on the asymptotic stability of positive equilibria have been shown in \cite{Mincheva.2004aa}.

In this work our focus is different from that of the literature cited above, although the asymptotic behavior of complex-balanced systems was part of our motivation (see Appendix \ref{sec:complexBal}). Namely, we are concerned with the convergence of a certain space-discretization scheme --the so-called method of lines-- for mass-action reaction-diffusion systems. We adopt the framework for convergence analysis introduced by Verwer \cite{Verwer.1984aa}, and concentrate on the proof-of-concept reaction
\begin{equation}\label{eq:reac}
A + B \overset{k_{-1}}{\underset{k_1}{\leftrightharpoons}} C
\end{equation}
within 1D space, while at the same time noting that our techniques are readily generalizable to other reaction-diffusion networks and to more than one space dimension. Indeed, it will be obvious how to extend our proofs to the multi-dimensional case; we only note that the proof of the comparison principle (the continuous and the discrete versions; see Section  \ref{sec:comparison}) imposes a limitation on the spatial dimension (should be at most five; see \cite{WCH} for details).

The Method of Lines (MOL) is not a mainstream numerical tool and the specialized literature is rather scarce. The method amounts to discretizing evolutionary PDE's in space only, so it produces a semi-discrete numerical scheme which consists of a system of ODE's (in the time variable). To prove convergence of the semi-discrete MOL scheme to the original PDE one needs to perform some more or less traditional analysis: it is necessary to show that the scheme is consistent with the continuous problem, and that the discretized version of the spatial differential operator retains sufficient dissipative properties in order to allow an application of Gronwall's Lemma to the error term.  As shown in \cite{Verwer.1984aa}, a uniform (in time) consistency estimate is sufficient to obtain convergence; however, the consistency estimate we proved is not uniform for small time, so we cannot directly employ the results in \cite{Verwer.1984aa} to prove convergence in our case. Instead, we prove all the required estimates ``from scratch'', then we use their exact quantitative form in order to conclude convergence.

For \eqref{eq:reac} we adopt the following paradigm: we envision splitting the spatial domain into $N$ equal subintervals in each of which we treat the concentrations of the three species as approximately constant. This, of course, is a fairly reasonable assumption for large $N$. We assume that a version of \eqref{eq:reac} takes place in each cell (or ``box'') $k$ (see Figure \ref{fig:boxes}), and the diffusion of any of the three species can be thought of as a reaction between adjacent replicas of the same species. The coefficients of these same-species reactions must be proportional to $N^2$  in order to get diffusion in the $N\rightarrow\infty$ limit (see also \cite{GSW} for an explanation of this scaling). The plan is to show that the standard reaction-diffusion system corresponding to \eqref{eq:reac} is obtained from these approximating reaction systems in the $N\rightarrow\infty$ limit.

The paper is organized as follows. The next section sets up the notation and preliminaries needed to state  our main result, Theorem  \ref{the0}.  In Section \ref{sec:comparison} we discuss comparison principles for solutions of the reaction-diffusion equation corresponding to \eqref{eq:reac} and for its space discretization, which we then use to prove consistency and boundedness of the logarithmic norm (in the spirit of \cite{Verwer.1984aa}, even though we had to make do with a nonuniform estimate). This completes the proof of
Theorem \ref{the0}, and it is done in Section \ref{sec:convergence}. The appendix collects a few technical results regarding the heat kernel and needed in the proof of Theorem \ref{the0}. MOL has an interesting interpretation in the context of mass-action reaction-diffusion systems, and particularly for complex-balanced networks. This is explored at the end of the appendix (Section \ref{sec:complexBal}), and discussed in connection with asymptotic results from literature and future directions of work.


\section{Main result}

Let $I:=(0,1)$. The primary concern of this work is the system of semi-linear parabolic partial differential equations
\begin{equation*}\label{E1}\tag{$E1$}
\begin{cases}
\partial_t a(t,x) = -k_1 a(t,x) b(t,x) + k_{-1}c(t,x) + k_A\,\ \partial_x^2 a(t,x) & \mbox{in $[0,T)\times I$}, \\
\partial_t b(t,x) = -k_1 a(t,x) b(t,x) + k_{-1}c(t,x) + k_B\,\ \partial_x^2 b(t,x)  & \mbox{in $[0,T)\times I$},\\
\partial_t c(t,x) = k_1 a(t,x) b(t,x) - k_{-1}c(t,x) + k_C\,\ \partial_x^2 c(t,x) & \mbox{in $[0,T)\times I$}
\end{cases}
\end{equation*}
together with the homogeneous Neumann boundary conditions
\begin{equation}\label{NC}
\begin{cases}
\partial_x a(t,0)=\partial_x a(t,1)=0,\\
\partial_x b(t,0)=\partial_x b(t,1)=0, \\
\partial_x c(t,0)=\partial_x c(t,1)=0.
\end{cases}
\end{equation}

 Here $k_1, k_{-1}$ and $T$ are positive constants and $k_A, k_B$ and $ k_C$ are the constant positive diffusion coefficients.

 The problem should be well-posed once appropriate initial conditions
$a(0,x)= a_0(x), \,\,\,\, b(0,x)= b_0(x) \,\,\,\,\mbox{and}\,\,\,\,
c(0,x)= c_0(x)$ are given.
 In the case of reaction-diffusion systems, there are two different aspects of existence to consider: local (in time) existence and global (in time) existence of solutions. The existence question is, in general, difficult to deal with. The well-posedness for  a general form of nonlinear parabolic system was obtained in \cite{CE}.
 In addition, they established existence and uniqueness for specific, three species systems when the diffusion coefficients are the same for all three species. In \cite{WCH} the authors established global existence and uniqueness of  solutions to \eqref{E1} with constants $k_1= k_{-1}= 1$ and distinct diffusion coefficients $k_A, k_B, k_C .$

  \subsection{Discretization by the Method of Lines (MOL)} We now discretize \eqref{E1} in space only: more precisely, we use the standard three-point stencil to approximate the second-order spatial derivatives.
Let $N\geq 2$ and divide the interval $[0,1]$ into $N$ subintervals of equal length $h:=1/N$, so that we have $N+1$ mesh points spaced by $h$ and numbered from $0$ to $N$.
The discretized problem is
\begin{equation}\label{E2}\tag{$E2$}
 \begin{cases}
\dot a_k^N(t) = -k_1 a_k^N(t) b_k^N(t) + k_{-1}c_k^N(t)+ (k_A/h^2)[a_{k-1}^N(t) - 2 a_k^N(t) +  a_{k+1}^N(t)] \\
\dot b_k^N(t) = -k_1 a_k^N(t) b_k^N(t) + k_{-1}c_k^N(t)+ (k_B/h^2)[b_{k-1}^N(t)- 2 b_k^N(t) +  b_{k+1}^N(t)] \\
\dot c_k^N(t) = k_1 a_k^N(t) b_k^N(t) - k_{-1}c_k^N(t)+ (k_C/h^2)[c_{k-1}^N(t)- 2 c_k^N(t) +  c_{k+1}^N(t)] \\
a^N_k(0) = N \int^{k/N}_{(k-1)/N} a_0(x)\,dx \\
b^N_k(0) = N \int^{k/N}_{(k-1)/N} b_0(x)\, dx\\
c^N_k(0) = N \int^{k/N}_{(k-1)/N} c_0(x)\, dx\\
a_0^N (t):= a_1^N(t),\  b_0^N (t):= b_1^N(t) \  \mbox{and} \  c_0^N (t):= c_1^N(t),\\
a_{N+1}^N(t):= a^N_N(t),\  b_{N+1}^N(t):= b^N_N(t) \  \mbox{and} \  c_{N+1}^N(t):= c^N_N(t)
\end{cases}
\end{equation}
for $k = 1,...,N$, and $( ^.  = \frac{d}{dt}).$
For the left endpoint $x=0$ ($k=0$) we use the forward difference approximation
$$\partial_xa(0,t) \approx \frac{a^N_{1}(t)- a^N_0(t)}{h}= 0,$$
so we assume $a^N_0(t)= a^N_1(t)$. For the right endpoint $x=1$ ($k=N$) we use the backward difference approximation $$\partial_xa(1,t) \approx \frac{a^N_{N}(t)- a^N_{N+1}(t)}{h}= 0,$$
which gives $a^N_{N+1}(t)= a^N_N(t)$. The same holds for $b$ and $c$. Let $$\vec {u}^N(t) := \left[ \vec{a}^N(t)^T,\vec{b}^N(t)^T, \vec{c}^N(t)^T\right]^T  \in \mathbb R^{3N}$$
denote a solution of \eqref{E2} with the column vector
$\vec{a}^N(t) = [ a_1^N(t), \dots, a_N^N(t)] ^T\in \mathbb R^{N},$ and similar definitions for $\vec{b}^N(t), \vec{c}^N(t)\in \mathbb R^{N}.$ (Note that $\vec {u}^N(t)$ is a column vector as well.)

  What we presented above is known as {\it Method of Lines} (MOL) \cite{Sadiku}; this nomenclature comes from the fact that we have reduced the original problem of finding a solution for \eqref{E1} at all points in the space-time rectangular domain $I\times[0,T]$ to the problem of finding a solution $\vec{u}^N$ on a finite number of lines in the space-time domain. By this method we store the concentrations at $N+1$ mesh points spaced by $h$ and numbered $0$ to $N$, and estimate the second derivatives of these concentrations at every point by using these values. The result of carrying out this procedure is a discretization of the system. The discretization is a set of ODEs \eqref{E2} which formally reduce to the original PDE \eqref{E1} in the $N\rightarrow\infty$ limit. Note that this method is also called semi-discretization because \eqref{E1} is discretized in space only.

\begin{figure}
%
%
%
\includegraphics[scale=.17]{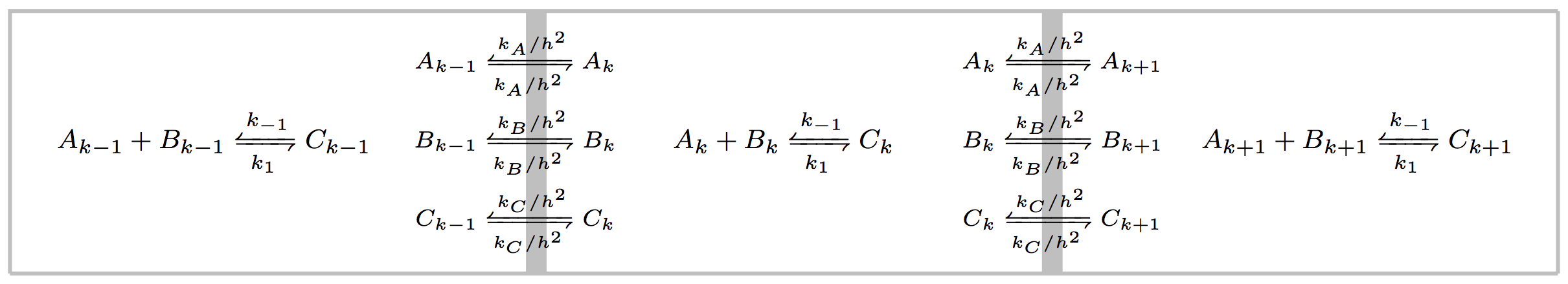}
\caption{Boxes $k-1$, $k$ and $k+1$}\label{fig:boxes}
\end{figure}

The setup of MOL described above is particularly intuitive for chemical networks. The space is divided into $N$ equal ``boxes" of homogeneous chemical compositions, with species transitions between adjacent boxes accounting for diffusion; see Figure \ref{fig:boxes}. The result is a reaction network with $3N$ species whose mass-action dynamics is given by (\ref{E2}). Clearly, the construction outlined here can be done starting from any reaction network, and it is discussed in Appendix \ref{sec:complexBal}.

We now define the following functions in a piecewise fashion. For $t>0$ let
$$a^N(t,x):=a^N_k(t),\ b^N(t,x):=b^N_k(t),\ c^N(t,x):=c^N_k(t)\ \mbox{ if } (k-1)/N\leq x<k/N.$$
Throughout the paper we denote by $\|\cdot\|_2$ the $L^2(0,1)$-norm, and by $|\cdot|$ and $\langle\cdot, \cdot \rangle$ the Euclidian norm and inner product in any ${\mathbb R}^n$.

We are now ready to present the main result of our paper:
 \begin{theorem}\label{the0}
 Assume that the initial data $a_0, b_0 $ and $c_0$ in $L^{\infty}(0,1)$ such that $a_0 \geq 0, b_0 \geq 0,$ and $c_0 \geq 0$ a.e. in $(0,1)$. Then the solution $\vec{u}^N $ of \eqref{E2} converges in $L^2(0,1)$ to the solution  $(\alpha,\beta,\gamma)$ of \eqref{E1}
in the sense
\begin{equation}\label{main-cv}
\lim_{N\rightarrow \infty}\bigg[\|a^N(t,\cdot)-\alpha(t,\cdot)\|_2+\|b^N(t,\cdot)-\beta(t,\cdot)\|_2 +\|c^N(t,\cdot)-\gamma(t,\cdot)\|_2\bigg] =0
\end{equation}
$\mbox{for all}\,\ t \in [0,T).$
\end{theorem}

\section{The Comparison Principle (Continuous and Discrete Problems)}\label{sec:comparison}

Chen, Li and Wright \cite{WCH} established a maximum principle for a version of \eqref{E1} on the whole real line and the constants $k_{-1}=k_1= 1, k_A > 0, k_B > 0, k_C > 0$. We adapt their proof to our case; the adaptation (Theorem \ref{the-bound-thm} below) is quite straightforward but we show its proof in some detail, mainly because the same proof will work for the discrete problem \eqref{E2} if one replaces the Heat Kernel by its discrete version.\\

 We change the variables by setting
$$u_1(t,x):=k_1 a(t,x), \,\ u_2(t,x):=k_1 b(t,x), \,\, u_3(t,x):=k_1 c(t,x),$$
so that \eqref{E1} becomes
\begin{equation}\label{P}
\begin{cases}
\partial_t u_1(t,x) = - u_1(t,x) u_2(t,x) +k_{-1} u_3(t,x) +k_A \,\ \partial_x^2 u_1(t,x), \\
\partial_t u_2(t,x)= -u_1(t,x) u_2(t,x) +k_{-1} u_3(t,x) +k_B \,\partial_x^2 u_2(t,x), \\
\partial_t u_3(t,x)= \,\    u_1(t,x) u_2(t,x) -k_{-1} u_3(t,x) +k_C \,\ \partial_x^2 u_3(t,x).
\end{cases}
\end{equation}
and the boundary conditions become
\begin{equation}\label{NC'}
\begin{cases}
\partial_x u_1(t,0)=\partial_x u_1(t,1)=0,\\
\partial_x u_2(t,0)=\partial_x u_2(t,1)=0, \\
\partial_x u_3(t,0)=\partial_x u_3(t,1)=0.
\end{cases}
\end{equation}
We also know that (see \cite{WCH} or \cite{Rothe.1984aa}) the solutions $u_1, u_2, u_3$ stay nonnegative if the initial data $u_{1,0}, u_{2,0}$ and $ u_{3,0}$ are nonnegative.

The following lemma, adapted from \cite{WCH}, will be used to prove Theorem \ref{the-bound-thm}.
\begin{lemma}\label{lemma1}
Assume that $u_{1,0}, u_{2,0}$ and $u_{3,0}$ are nonnegative. Then there exists a constant $C$ such that $$\int_0^1  u_i(t,x)  \ dx \leq C \quad\mbox{for all $t$ and for $i= 1,2,3.$}$$
\end{lemma}
\proof  Equations \eqref{P} yield
\begin{equation*}
\partial_t [u_1(t,x) + u_2(t,x) + 2u_3(t,x)] = k_A \,\partial_x^2 u_1(t,x)+ k_B \,\partial_x^2 u_2(t,x) + 2 k_C \,\partial_x^2 u_3(t,x),
\end{equation*}
and therefore
\begin{align*}
\partial_t &\int ^1_0 [u_1(t,x) + u_2(t,x) + 2u_3(t,x)]\,dx =\\
&= k_A \, \int ^1_0 \partial_x^2 u_1(t,x) \,dx + k_B \, \int ^1_0 \partial_x^2 u_2(t,x)\,dx + 2 k_C \, \int ^1_0\partial_x^2 u_3(t,x)\,dx,
\end{align*}
or, using the boundary conditions \eqref{NC'}:
\begin{equation}
\partial_t \int ^1_0 [u_1(t,x) + u_2(t,x) + 2u_3(t,x)]\,dx = 0.
\end{equation}
It follows that
\begin{equation*}
\int ^1_0 [u_1(t,x) + u_2(t,x) + 2u_3(t,x)]\,dx = \int ^1_0 [u_{1,0}(x) + u_{2,0}(x) + 2u_{3,0}(x)]\,dx = C
\end{equation*}
is constant. Since $u_1, u_2, u_3$ are nonnegative \cite{WCH}, we end up with
$$\int_0^1  u_i(t,x) \ dx \leq C $$ for all $t$ and $i= 1,2,3.$
\endproof

\begin{theorem}\label{the-bound-thm}
If the initial data $a_0,\,b_0,\,c_0$ of the system \eqref{P} belong to $L^\infty(0,1)$, then the solution of \eqref{P} is uniformly bounded for all time.
\end{theorem}
\proof Let $H_i(t,x,y)$ denote the Neumann Heat Kernel of the linear parabolic equation
$$\partial_t v = \mathcal{L}_i v,$$ where
$\mathcal{L}_1=k_A\partial_x^2$, $\mathcal{L}_2=k_B\partial_x^2$ and $\mathcal{L}_3=k_C\partial_x^2.$
Note that $H_1(t,x,y)=H(k_At,x,y)$, $H_2(t,x,y)=H(k_Bt,x,y)$, and $H_3(t,x,y)=H(k_Ct,x,y),$ where $H$ is the Neumann Heat Kernel defined in Appendix.

First, from the nonnegativity of the solutions, we have
$$\partial_t u_1- k_A \,\partial_x^2 u_1= - u_1 u_2+ k_{-1} u_3,$$
which implies
\begin{equation}
 \partial_t u_1- k_A\,\partial_x^2 u_1 \leq  k_{-1} u_3.
\end{equation}
For each fixed $T>0$, we compare $u_1(T+\cdot,\cdot)$ with the solution of the linear equation
\begin{equation}\label{eq31}
\begin{cases}
\partial_t v -\mathcal{L}_1 v = k_{-1} u_3(T+t,x)\\
v(0,x)=u_1(T,x).
\end{cases}
\end{equation}
We know that the solution of \eqref{eq31} is
\begin{align*}
v(T+t,x) &= \int_I H(k_At,x,y)\,u_1(y,T)\,dy +\\
&+\int_0^t \int_I H(k_A(t-s),x,y)\,k_{-1}\,u_3(T+s,y)\,dy \,ds.
\end{align*}
With $I_1$ and $I_2$ denoting the first and second integral terms above, it now follows that
$$u_1(T+t,x) \leq I_1 + I_2.$$
When $t$ is bounded away from $0$, $H(t,x,y)$ is bounded in a pointwise sense. In the following, we assume that
$\delta/2 \leq t \leq \delta \mbox{ for some } \delta > 0.$
According to \eqref{BH} (Appendix), the Neumann Heat Kernel satisfies the bounds
\begin{equation}\nonumber
0\leq H(t,x,y)\leq  1+ 2 f(4 t)\leq 1+2f(2\delta)\mbox{ for all }x,\ y\in I,\ t\geq\frac{\delta}{2},
\end{equation}
where $f$ is defined in Appendix \eqref{special_fun}. Therefore, the integral $I_1$ can be easily bounded as
\begin{equation}\nonumber
I_1 = \int_I H_1(t,x,y)\,u_1(T,y)\,dy  \leq [ 1+2 f(2 k_A\delta)] \int_I   u_1(T,y)\,dy \leq C [ 1+2 f(2k_A \delta)],
\end{equation}
where, by Lemma \ref{lemma1}, $C$ is a constant which is independent of $T$.
As for the integral $I_2$, we can rewrite it as
\begin{align*}
I_2 &= \int_0^t \int_I H_1(t-s,x,y)\,k_{-1}\,u_3(T+s,y)\,dy \,ds\\
&=  \int_0^{t-\frac{\delta}{2}} \int_I H_1(t-s,x,y)\,k_{-1}\,u_3(T+s,y)\,dy \,ds\ + \\
&+\int_{t-\frac{\delta}{2}}^t \int_I H_1(t-s,x,y)\,k_{-1}\,u_3(T+s,y)\,dy \,ds\\ &{=: } \hat{I_2} + \tilde{I_2}.
\end{align*}
Next we estimate
\begin{align*}
\hat{I_2} &= \int_0^{t-\frac{\delta}{2}} \int_I H_1(t-s,x,y)\,k_{-1}\, u_3(T+s,y)\,dy \,ds\\ &\leq  \int_0^{t-\frac{\delta}{2}} ( 1+2 f(4k_A(t-s))) k_{-1} \int_I u_3(T+s,y)\,dy\,ds\\
&\leq   k_{-1} C  \int_0^{t-\frac{\delta}{2}} ( 1+2 f(4 k_A(t-s)))\,ds\\ & \leq \frac{k_{-1} C\delta}{2}+\frac{k_{-1}C}{2k_A}\int_{2 k_A\delta}^{4k_A\delta}f(\tau)d\tau=O(\delta).
\end{align*}
For the integral $\tilde{I_2}$, we use H\"older's inequality to get, for any $\infty>p,\ q>1$ such that $1/p+1/q=1$,
\begin{align*}
\tilde{I_2} &= \int_{t-\frac{\delta}{2}}^t \int_I H_1(t-s,x,y)\,k_{-1}\,u_3(T+s,y)\,dy \,ds \\
&\leq \int_{t-\frac{\delta}{2}}^t \left| \int_I H_1(t-s,x,y)\,k_{-1}\,u_3(T+s,y)\,dy\right| \,ds\\
&\leq k_{-1} \int_{t-\frac{\delta}{2}}^t \left( \int_I \big| H_1(t-s,x,y)\big|^q\,dy\right)^{1/q}\, \left(  \int_I \left| \,u_3(T+s,y)\right|^p\,dy\right)^{1/p} \,ds\\
&\leq k_{-1} \int_{t-\frac{\delta}{2}}^t \left( \int_I \left| 1+2 f(4k_A(t-s)) \right|^q\,dy\right)^{1/q}\, \left( \int_I \left|\,u_3(T+s,y)\right|\,\left| \,u_3(T+s,y)\right|^{p-1}\,dy \right)^{1/p} \,ds\\
&\leq k_{-1}  \int_{t-\frac{\delta}{2}}^t  \,\left| 1+2 f(4k_A(t-s)) \right|  \left( \left\| \,u_3(T+s,\cdot)\right\|^{p-1}_{ \infty} \int_I \left|\,u_3(T+s,y)\right|\,dy \right)^{1/p} \,ds\\
&\leq  C k_{-1} \smash{\displaystyle\max_{t-\frac{\delta}{2} \leq s \leq t}} \left\| u_3(T+s,\cdot)\right\|^{\frac{p-1}{p}}_{\infty}
\int_{t-\frac{\delta}{2}}^t \left[1+2 f\left(4k_A\left(t-s\right)\right)\right]  \,ds\\
&\leq  C k_{-1}  \left( \frac{ \delta}{2}+\frac{1}{2k_A}\int_0^{2 k_A\delta}f(\tau)d\tau\right) \smash{\displaystyle\max_{t-\frac{\delta}{2} \leq s \leq t}} \left\| u_3(T+s,\cdot)\right\|^{\frac{p-1}{p}}_{\infty}.
\end{align*}
Therefore, we get
 \begin{equation*}
u_1(T+t,x)\leq C + C \smash{\displaystyle\max_{t-\frac{\delta}{2} \leq s \leq t}} \lVert u_3(T+s,\cdot)\rVert^{\frac{p-1}{p}}_{\infty},
\end{equation*}
where the constant $C_1$  depends on $\delta$ (in fact, due to the integrability of $f$, $C_1$ tends to zero as $\delta$ tends to zero). Since $\frac{\delta}{2} \leq t \leq \delta$, we deduce
\begin{equation}\label{eq5}
\smash{\displaystyle\max_{T+ \frac{\delta}{2} \leq s \leq T+ \delta}} \lVert u_1(s,\cdot)\rVert_{\infty} \leq C + C \smash{\displaystyle\max_{T \leq s \leq T+ \delta}} \lVert u_3(s,\cdot)\rVert^{\frac{p-1}{p}}_{\infty}.
\end{equation}
Applying the same argument to the equation for $u_2$ from the system \eqref{P}, we also have
\begin{equation}\label{eq6}
\smash{\displaystyle\max_{T+ \frac{\delta}{2} \leq s \leq T+ \delta}} \lVert u_2(s,\cdot)\rVert_{\infty} \leq C + C \smash{\displaystyle\max_{T \leq s \leq T+ \delta}} \lVert u_3(s,\cdot)\rVert^{\frac{p-1}{p}}_{ \infty}.
\end{equation}
 Finally, from the third equation of the system \eqref{P}, we have
$$\partial_t u_3- k_C \,\partial_x^2 u_3=  u_1 u_2- k_{-1} u_3,$$
which implies
\begin{equation}\nonumber
 \partial_t u_3- k_C  \,\partial_x^2 u_3 \leq  u_1 u_2.
\end{equation}
Just as before, one gets
\begin{equation}\nonumber
u_3(T+t,x)\leq C + C \smash{\displaystyle\max_{t-\frac{\delta}{2} \leq s \leq t}} \lVert u_1(T+s,\cdot)\rVert_{ \infty} \lVert u_2(T+s,\cdot)\rVert_{ \infty}^{(p-1)/p},
\end{equation}
which yields
\begin{equation}\label{eq7}
\smash{\displaystyle\max_{T+ \frac{\delta}{2} \leq s \leq T+ \delta}} \lVert u_3(s,\cdot)\rVert_{\infty}\leq C + C \smash{\displaystyle\max_{T \leq s \leq T+ \delta}} \lVert u_1(s,\cdot)\rVert_{ \infty} \lVert u_2(s,\cdot)\rVert_{ \infty}^{(p-1)/p}.
\end{equation}
As in \cite{WCH}, we can use  (\ref{eq5})--(\ref{eq7}) to get
\begin{equation}\label{eq8}
\smash{\displaystyle\max_{T+ \frac{\delta}{2} \leq s \leq T+ \delta}} \lVert u_3(s,\cdot)\rVert_{\infty}\leq C + C\smash{\displaystyle\max_{T- \frac{\delta}{2} \leq s \leq T+ \delta}} \lVert u_3(s,\cdot)\rVert _{ \infty}^{\omega},
\end{equation}
where
$$\omega:= \frac{p-1}{p} + \bigg(\frac{p-1}{p}\bigg)^2.$$ We choose
$$1<p<\frac{2}{3-\sqrt{5}}, \mbox{ so that }0<\omega<1.$$
Once more, as in \cite{WCH}, we infer
\begin{equation}\label{eq9}
\smash{\displaystyle\max_{T+ \frac{\delta}{2} \leq s \leq T+ \delta}} \lVert u_3(s,\cdot)\rVert_{\infty}\leq C + C \Big( \smash{\displaystyle\max_{T- \frac{\delta}{2} \leq s \leq T+ \frac{\delta}{2}}} \lVert u_3(s,\cdot)\rVert _{ \infty}^{\omega} +   \smash{\displaystyle\max_{T+ \frac{\delta}{2} \leq s \leq T+ \delta}} \lVert u_3(s,\cdot)\rVert_{ \infty}^{\omega}\Big),
\end{equation}
which implies
\begin{equation}\label{eq10}
\smash{\displaystyle\max_{T+ \frac{\delta}{2} \leq s \leq T+ \delta}} \lVert u_3(s,\cdot)\rVert_{ \infty}\leq C + C_2 \smash{\displaystyle\max_{T- \frac{\delta}{2} \leq s \leq T+ \frac{\delta}{2}}} \lVert u_3(s,\cdot)\rVert _{ \infty}^{\omega}.
\end{equation}
Denote
\begin{equation*}
 M(t):=  \smash{\displaystyle\max_{t-\frac{\delta}{2} \leq s \leq t+ \frac{\delta}{2}}} \lVert u_3(s,\cdot)\rVert_{ \infty}^\omega.
 \end{equation*}
Then, one can derive (see \cite{WCH}) from (\ref{eq10}) that
$$ M(t)\leq C_2 + C_3 \left[ M\left(t- \frac{\delta}{2}\right)\right]^\omega\mbox{ for all }\frac{\delta}{2}\leq t\leq\delta.$$
Since $0<\omega<1$, we deduce $u_3$ is bounded for all time if $M(\frac{\delta}{2})< \infty$, and thus, $u_1$ and $u_2$ are also bounded for all time. This completes the proof of the theorem.
\endproof

Note that the proof of the Theorem \ref {the-bound-thm} goes through if we replace $H$ by $H^N$ (i.e. the discrete Neumann Heat Kernel on $I$). Indeed, as shown in Appendix, Subsection \ref{special_fun}, $H^N$ has all the desired properties.

\begin{theorem}\label{the2}
Assume that the initial data $a_0,\,b_0,\,c_0$ in the system \eqref{E2} are essentially bounded by a positive constant $M$. Then there exists a finite constant $\tilde{M}$ such that for all $N$, $k=1,2, ..., N$ and $t\in (0, \infty)$ we have
$$a^N_k(t), b^N_k(t), c^N_k(t)\in [0,\tilde M]$$
\end{theorem}
\proof We can apply the same proof as above, by replacing the Heat Kernel $H$ by $H^N$ and using \eqref{eqa^N}, \eqref{eqb^N} and \eqref{eqc^N}.
\endproof

\section{Convergence}\label{sec:convergence}

 In this section we will prove the main result. We first need to check the {\it consistency} of the MOL scheme when applied to our system.

We solve \eqref{E1} with initial data $a_0,\, b_0,\, c_0$, and denote by $(\alpha,\beta,\gamma)$ a solution. Let $N\geq 2$ be integer. For each $t\geq0$, we define the column vectors $\vec\alpha^N(t),\, \vec\beta^N(t),\, \vec\gamma^N(t)$ whose components are
\begin{align*}
\alpha^N_k(t):=\alpha(t,(k-1)/N)\mbox{ for all }k=1,...,N,\\
\beta^N_k(t):=\beta(t,(k-1)/N)\mbox{ for all }k=1,...,N,\\
\gamma^N_k(t):=\gamma(t,(k-1)/N)\mbox{ for all }k=1,...,N.
\end{align*}
For $t>0$ and  $(k-1)/N\leq x<k/N$ define
\begin{equation}\label{alphaN}\alpha^N(t,x):=\alpha^N_k(t),\ \beta^N(t,x):=\beta^N_k(t),\ \gamma^N(t,x):= \gamma^N_k(t)
\end{equation} and let
\begin{equation}\label{solution-vector}
\vec{v}^N(t)=\left[\vec\alpha^N(t)^T,\vec \beta^N(t)^T,\vec\gamma^N(t)^T\right]^T \in \mathbb{R}^{3N}.
\end{equation}

 Also, let us denote the discrete Laplacian matrix with Neumann boundary condition on $I$ by
\begin{equation}\label{Delta^N}
\Delta^N :=  N^2 \begin{bmatrix}
-1 & 1 & 0 & 0 & 0 & \ldots & 0 \\
1 & -2 & 1 & 0 & 0 & \ldots & 0\\
0 & 1 & -2 & 1 & 0 & \ldots & 0 \\
\vdots & \vdots & \vdots & \vdots & \vdots & \vdots & \vdots \\
0 & 0 & \ldots & 0 & 1 & -2 & 1 \\
0 & 0 & \ldots & 0  & 0 & 1 & -1
\end{bmatrix} _ {N \times N}.
\end{equation}

Next, consider the vector field $\widetilde{F}^N: \mathbb{R}^{3N} \rightarrow \mathbb{R}^{3N}$ given by
\[
\widetilde{F}^N = \begin{bmatrix} F^N \\ F^N \\ -F^N \end{bmatrix},
\]
where $F^N:\mathbb{R}^{3N} \rightarrow \mathbb{R}^N,$
$$F^N ([\vec a^T, \vec b^T, \vec c^T]^T) = \begin{bmatrix}
-k_1 a_1 b_1 + k_{-1}c_1 \\
\vdots \\
-k_1 a_N b_N + k_{-1}c_N \\
 \end{bmatrix}
$$
for  generic column vectors $\vec a, \vec b, \vec c \in {\mathbb R}^{N}.$

Let $\tilde \Delta^N$ be the $(3N)\times(3N)$ block-diagonal matrix whose diagonal blocks are the matrices $k_A\Delta^N,\, k_B\Delta^N$ and $k_C\Delta^N$. Note that the discrete system \eqref{E2} can now be written as
\begin{equation}\label{eq1}
\frac{d}{dt}\vec{u}^N(t)=\widetilde{F}^N\left(\vec{u}^N(t)\right)+ \tilde\Delta^N\vec{u}^N(t)
\end{equation}
and consists of three coupled systems
\begin{align*}\label{eq2}
\dot{\vec{a}}^N(t)&= {F}^N\left(\vec{u}^N(t)\right)+k_A \Delta^N \vec{a}^N(t),\\
\dot{\vec{b}}^N(t)&= {F}^N\left(\vec{u}^N(t)\right)+k_B\Delta^N \vec{b}^N(t),\\
\dot{\vec{c}}^N(t)&= - {F}^N\left(\vec{u}^N(t)\right)+k_C\Delta^N \vec{c}^N(t).
\end{align*}

%

\subsection{A consistency estimate}
We begin by proving an estimate on the {\it space truncation error}; this is called a {\it consistency estimate}.
  This is the error obtained by ``plugging'' the solution to the continuous problem \eqref{E1} into the approximating discrete scheme.
 In order to do that, note that we can write a system of equations for $\vec{v}$ (defined in \eqref{solution-vector}) in the form
\begin{equation}\label{eqv}
\frac{d}{dt}\vec{v}^N(t)=\widetilde{F}^N\left(\vec{v}^N(t)\right)+\tilde\Delta^N\vec{v}^N(t)+\vec \varepsilon^N(t).
\end{equation}
Here $\vec \varepsilon^N(t):=\left[\vec\varepsilon^{N,\alpha}(t),\vec\varepsilon^{N,\beta}(t),\vec\varepsilon^{N,\gamma}(t)\right]^T\in \mathbb{R}^{3N}$, where
$\vec\varepsilon^{N,\alpha}\in \mathbb{R}^{N}$ has components
\[
\quad \varepsilon^{N,\alpha}_k(t):=k_A\partial^2_x\alpha \left(t,(k-1)/N \right) - k_A N^2 \left[ \alpha^N_{k-1}(t) - 2 \alpha^N_k(t) + \alpha_{k+1}^N(t) \right],
\]
and $\vec\varepsilon^{N,\beta},\ \vec\varepsilon^{N,\gamma}$ are defined similarly.
It is readily seen that
\begin{equation}\label{eqe2}
\big|\varepsilon^{N,\alpha}_k(t)\big| \leq \frac{1}{N} \left\|\partial^3_x\alpha(t,\cdot)\right\|_{L^\infty(0,1)},
\end{equation}
and similarly for $\beta$, $\gamma$; this implies that  a consistency estimate boils down to bounding the third spatial derivates, which we pursue next:

\begin{theorem}\label{thm-deriv-bd}
 Let $a_0,\,b_0,\,c_0  \in L^\infty(0,1)$ and $a_0 \geq 0,\,b_0 \geq 0,\,c_0 \geq 0$ a.e. in $I$. Let $k_A,\,k_B,\,k_C \geq 0$ and $k_{-1},\,k_1 >0$. Consider the solution $(\alpha,\beta,\gamma)$ for the system \eqref{E1}
with initial data $a_0,\,b_0,\,c_0$. Fix $0<T< \infty$. Then
 for any integer $j\geq 1$ and any $0<\delta <T$ the derivatives $\partial_x^j \alpha,\,\partial_x^j \beta,\, \partial_x^j \gamma$ belong to $L^\infty((\delta,T)\times I)$, with norm upper bounds depending only on $j$, $\delta,\ T,\ M,\ k_{-1},\ k_1$.
\end{theorem}
\proof
By Duhamel's Principle, we have
\begin{align}\label{D-P}
\alpha(t,x) = &\int_I H\left(k_A(t-\delta),x,y\right)\, \alpha(\delta,y)\,dy \\
\nonumber &+ \int_\delta^t \int_I H \left(k_A(t-s),x,y\right) \Big[k_{-1}\gamma(s,y) -k_1 \alpha(s,y)\,\beta(s,y)\Big]\,dy\,ds
\end{align}
for $t>\delta\geq 0$. It is easy to see that
\[
\int_I \Big|\partial_y H \left(k_At,x,y\right)\Big|\,dy \leq C(k_A,\delta,T) < \infty
\]
for all $\delta \leq t <T$,uniformly in $x\in I$.

Proposition \ref{deriv-integr} also guarantees, in light of the property $(5')$ of $H$ (see Appendix \ref{sec:heatKernels}), that, for all $t\in(0,T]$, the kernel $H(t,x,y)$ satisfies (uniformly in $x\in I$)
\begin{equation}\nonumber
\int_\delta^t \int_I \Big|\partial_x H\left(k_A(t-s),x,y\right)\Big|\,dy\,ds  \leq C(k_A,\delta,T) < \infty.
\end{equation}
Therefore, we can differentiate under the integral in \eqref{D-P} to see that, if $2\delta\leq t\leq T$, then $\alpha(t,\cdot)$ is differentiable on $(0,1)$, and
\begin{align*}
\partial_x \alpha(t,x) = &\int_I \partial_x H \left( k_A(t-\delta),x,y\right)\,\alpha(\delta,y)\,dy \\
\nonumber &+ \int_\delta^t \int_I \partial_x H \left(k_A(t-s),x,y\right) \Big[k_{-1}\gamma(s,y) - k_1 \alpha(s,y)\,\beta(s,y)\Big]\,dy\,ds.
\end{align*}

Now, using property (5') of $H$ again, and  and replacing $\partial_x H$ by $-\partial_y H_D$, we get
\begin{align}\label{d_x}
\partial_x \alpha(t,x) = &- \int_I \partial_y H_D\left(k_A(t-\delta),x,y\right)\,\alpha(\delta,y)\,dy \\
\nonumber &- \int_\delta^t \int_I \partial_y H_D\left(k_A(t-s),x,y\right)\Big[k_{-1}\gamma(s,y) -k_1 \alpha(s,y)\,\beta(s,y)\Big]\,dy\,ds.
\end{align}
Therefore, $\alpha(t,\cdot)$ is differentiable on $(0,1)$. The Neumann boundary conditions also give that $\partial_x \alpha(t,0) = \partial_x \alpha(t,1) =0$, so $\alpha(t,\cdot)$ is differentiable on $I$, with zero slopes at boundary.

Obviously, $\beta$ and $\gamma$ enjoy the same regularity. Thus, we can integrate by parts (in space) \eqref{d_x} to get
\begin{align}\label{d_x(a)}
\partial_x \alpha(t,x) &= \int_I H_D\left(k_A(t-\delta),x,y\right)\,\partial_y \alpha(\delta,y)\,dy \\ \nonumber
&+ \int_\delta^t \int_I H_D\left(k_A(t-s),x,y\right)\Big[k_{-1}\partial_y \gamma(s,y)-k_1 \alpha(s,y) \partial_y \beta(s,y)-\\
&-k_1 \beta(s,y) \partial_y \alpha(s,y)\Big]\,dy\,ds,
\end{align}
where we used that $H_D(t,x,0) = H_D(t,x,1) =0$ for all $t>0,\, x\in I$.

Let $u(t,x):= \partial_x \alpha(t,x),\, v(t,x):=\partial_x \beta(t,x),\, w(t,x):=\partial_x \gamma(t,x)$. From \eqref{d_x(a)}, we get
\begin{align*}
\left| u(t,x)\right| &\leq \left| \int_I H_D(k_A(t-\delta),x,y)\,u(\delta ,y)\,dy\right| +\\
&+\int_\delta^t \Big\{k_{-1} \left| \int_I H_D(k_A(t-s),x,y) \, w(s,y) \, dy \right|+\\
&+ k_1 M \left| \int_I H_D(k_A(t-s),x,y) \, v(s,y)\,dy \right| +\\
&+ k_1 M \left| \int_I H_D(k_A(t-s),x,y) \, u(s,y)\, dy \right| \Big\}\, ds.
\end{align*}

Since for any $f\in L^\infty(I)$, $g(t,x):=\int_IH_D(k_At,x,y)f(y)dy$ satisfies $g(0,\cdot)=f$ and solves $\partial_t g-k_A\partial^2_xg=0$ with Dirichlet boundary conditions, we conclude that each term whose absolute value is taken in the right hand side of the above inequality is the solution of the Dirichlet problem originating from the indicated function and evaluated at a later time; by property (2) of $H_D$ (Appendix \ref{sec:heatKernels}) we conclude
\begin{align*}
\left| u(t,x)\right| \leq \| u(\delta,\cdot)\|_\infty + \int_\delta^t \Big\{ k_1 M \Big(\|u(s,\cdot)\|_\infty + \|v(s,\cdot)\|_\infty\Big) + k_{-1} \| w(s,\cdot)\|_\infty \Big\}\,ds.
\end{align*}
Let $\lambda(k_{-1},k_1,M):=\max\{k_1 M, k_{-1}\}$. Then, for all $x\in I$ we have
\begin{align*}
| u(t,x)| \leq \| u(\delta,\cdot)\|_\infty + \lambda(k_{-1},k_1,M) \int_\delta^t \Big(\|u(s,\cdot)\|_\infty + \|v(s,\cdot)\|_\infty + \| w(s,\cdot)\|_\infty \Big)\,ds.
\end{align*}

Likewise, we get
\begin{align*}
| v(t,x)| \leq \| v(\delta,\cdot)\|_\infty + \lambda(k_{-1},k_1,M) \int_\delta^t \Big(\|u(s,\cdot)\|_\infty + \|v(s,\cdot)\|_\infty + \| w(s,\cdot)\|_\infty \Big)\,ds,
\end{align*}
and
\begin{align*}
|w(t,x)|  \leq \| w(\delta,\cdot)\|_\infty + \lambda(k_{-1},k_1,M) \int_\delta^t \Big(\|u(s,\cdot)\|_\infty + \|v(s,\cdot)\|_\infty + \| w(s,\cdot)\|_\infty \Big)\,ds.
\end{align*}
By addition and an application of Gronwall's Lemma, we get
\begin{equation}\label{B1}
\left\| u(t,\cdot)\right\|_\infty  + \| v(t,\cdot)\|_\infty + \| w(t,\cdot)\|_\infty \leq C_1(k_{-1},k_1,\delta,T,M) \quad\mbox{for all } t\in [2\delta,T].
\end{equation}

We now return to \eqref{d_x(a)} and, using now that the Dirichlet Kernel satisfies (see Appendix \ref{sec:heatKernels})
\begin{align*}
\int_I \left| \partial_x H_D(t,x,y)\right|\,dy <\infty \quad\mbox{ for }t\geq 2\delta, \ \mbox{ uniformly in } x\in I,
\end{align*}
and
\begin{align*}
 \int_\delta^t \int_I \left| \partial_x H_D(t-s,x,y)\right|\,dy\,ds\leq C(k_A,\delta,T) < \infty,
\end{align*}
we conclude that we can differentiate again with respect to $x$ under the integral signs. Therefore, we have
\begin{align*}
\partial_x^2 \alpha(t,x) &= \int_I \partial_x H_D\left(k_A(t-\delta),x,y\right)\,u(\delta,y)\,dy + \\
&+ \int_\delta^t \int_I  \partial_x H_D\left(k_A(t-s),x,y\right)\Big[ k_{-1} w(s,y) -\\
&-k_1 \beta(s,y) u(s,y) -k_1 \alpha(s,y) v(s,y)\Big]\,dy\,ds.
\end{align*}
Next we replace $\partial_x H_D$ by $-\partial_y H$ (property (5) of $H$ in Section \ref{sec:heatKernels}), and note that even though $H(t,x,\cdot)$ is not equal to zero at $y=0,\,1$, we can still integrate by parts and get rid of the boundary terms because $u(s,\cdot),\,v(s,\cdot),\,w(s,\cdot)$ are all zero at $y=0,\,1$ for all $s\in [\delta, t]$. Therefore, we get
 \begin{align}\label{d_xx(a)}
\partial_x^2 \alpha(t,x) &= \int_I H(k_A(t-\delta),x, y)\partial_y u(\delta,y)\,dy +\\
&+ \int_\delta^t \int_I  H\left(k_A(t-s),x,y \right) \Big[ k_{-1} \partial_y w(s,y) -2 k_1 u(s,y) v(s,y) -\\ \nonumber
&-k_1 \beta(s,y) \partial_y u(s,y) -k_1 \alpha(s,y) \partial_y v(s,y) \Big]\,dy\,ds. \nonumber
\end{align}
Property $(2')$ of $H$ in Section \ref{sec:heatKernels} implies

\begin{align*}
\left\| \partial_x^2 \alpha(t,\cdot)\right\|_\infty &= \sup_{x\in I} \left| \partial_x^2 \alpha(t,x)\right| \leq  \left\| \partial_y^2 \alpha(\delta,\cdot)\right\|_\infty
+\int_\delta^t \Big[ k_{-1}\left\| \partial_y^2 \gamma(s,\cdot)\right\|_\infty + \\
&+2 k_1 \left\| u(s,\cdot)\right\|_\infty \left\| v(s,\cdot)\right\|_\infty  + k_1 M \left\| \partial_y^2 \alpha(s,\cdot)\right\|_\infty +\\
&+ k_1 M \left\| \partial_y^2 \beta(s,\cdot)\right\|_\infty \Big]\,ds.
\end{align*}
We use \eqref{B1} to bound the term $\left\| u(s,\cdot)\right\|_\infty \left\| v(s,\cdot)\right\|_\infty$, write the corresponding inequalities for the $\beta$ and $\gamma$ terms, add them up and use Gronwall's Lemma again to get a bound $C_2(k_{-1},k_1,\delta,T,M)< \infty$ on $\left\| \partial_x^2 \alpha(t,\cdot)\right\|_\infty + \left\| \partial_x^2 \beta(t,\cdot)\right\|_\infty + \left\| \partial_x^2 \gamma(t,\cdot)\right\|_\infty$ for $t\in [2\delta,T]$.

From \eqref{d_xx(a)}, we differentiate again in $x$ to get (after using the property $(5')$ for $H$ yet again)
\begin{align*}
\partial_x^3 \alpha(t,x) &= - \int_I \partial_y H_D\big(k_A(t-\delta),x,y\big)\partial_y^2 \alpha(\delta,y)\,dy - \\
&-\int_\delta^t \int_I  \partial_y H_D\left(k_A(t-s),x,y\right)\Big[ k_{-1} \partial_y^2 \gamma(s,y)
-2 k_1 \partial_y \alpha(s,y) \partial_y \beta(s,y) - \\
&-k_1 \beta(s,y) \partial_y^2 \alpha(s,y) - k_1 \alpha(s,y) \partial_y^2 \beta(s,y)\Big]\,dy\,ds.
\end{align*}
This time we deal with the Dirichlet Kernel once more, so even if $\partial_y^2 \alpha(\delta,y)$ and the likes do not vanish at $y=0,\,1,\, H_D(t,x,\cdot)$ does for all $t>0$ and all  $x\in I$. Therefore, we can once more integrate by parts to get
\begin{align*}
\partial_x^3 \alpha(t,x) &=  \int_I  H_D\left(k_A(t-\delta),x,y\right)\partial_y^3 \alpha(\delta,y)\,dy + \\
&+\int_\delta^t \int_I  H_D\left(k_A(t-s),x,y\right)\Big[ k_{-1} \partial_y^3 \gamma(s,y) - 3 k_1 \partial_y^2 \alpha(s,y) \partial_y \beta(s,y) - \\
&-3 k_1 \partial_y \alpha(s,y) \partial_y^2 \beta(s,y) - k_1 \beta(s,y) \partial_y^3 \alpha(s,y) - k_1 \alpha(s,y) \partial_y^3 \beta(s,y)\Big]\,dy\,ds,
\end{align*}
which implies
\begin{align*}
\left| \partial_x^3 \alpha(t,x) \right| &\leq  \left\| \partial_y^3 \alpha(\delta,\cdot)\right\|_\infty + \int_\delta^t \Big\{ k_{-1}\left\| \partial_y^3 \gamma(s,\cdot)\right\|_\infty + \widetilde{C}(k_{-1},k_1,M,T,\delta) \\
&+  k_1 M \left\|  \partial_y ^3 \alpha(s,\cdot)\right\|_\infty + k_1 M \left\| \partial_y^3 \beta(s,\cdot)\right\|_\infty  \Big\}\,ds.
\end{align*}
Again, by Gronwall's Lemma, we get
$$\left\| \partial_x^3 \alpha(t,\cdot)\right\|_\infty + \left\| \partial_x^3 \beta(t,\cdot)\right\|_\infty + \left\| \partial_x^3 \gamma(t,\cdot)\right\|_\infty \leq C_3(k_{-1},k_1,M,T,\delta) <\infty.$$

The procedure can be continued to get bounds of the type
\begin{equation}\label{bdd}
\left\| \partial_x^{j} \alpha(t,\cdot)\right\|_\infty + \left\| \partial_x^{j} \beta(t,\cdot)\right\|_\infty + \left\| \partial_x^{j} \gamma(t,\cdot)\right\|_\infty \leq C(j,k_{-1},k_1,M,T,\delta) <\infty
\end{equation}
for all orders of differentiation $j\geq 1$. \endproof

\begin{remark}
The regularity assumed on initial data in Theorem \ref{thm-deriv-bd} prevents us from obtaining uniform bounds as $\delta\to 0$. Thus, Theorem 3.1 in \cite{Verwer.1984aa} cannot directly be applied here to yield a uniform (in time) consistency estimate.
\end{remark}

Since $k_1, k_{-1}$, and $M$ are fixed here, the bound in Theorem \ref{thm-deriv-bd} for third order derivatives only depends on $\delta$ and $T$. Denoting this quantity by $C(\delta,T)$, the consistency estimate now follows:

\begin{theorem} For any $0<\delta<T<\infty$ there exists a real constant $C(\delta,T)$ such that
\begin{equation}\label{eqeN1}
|\vec \varepsilon^N(t)|^2 \leq \frac{1}{N} C(\delta,T)\mbox{ for all integers }N\geq 2\mbox{ and all }t\in[\delta,T].
\end{equation}
\end{theorem}

\subsection{Proof of Theorem \ref{the0}}
Fix $T>0$ and $a_0 , b_0, c_0 \in L^\infty (0,1)$.
%
%
%
Let us begin by noticing that \eqref{main-cv} hold for $t=0$ (see Appendix, Subsection \ref{a0Ncv}). For $t\in(0,T)$ the proof is presented in three steps: first we prove that
$$\lim_{N\rightarrow\infty}\bigg[\|\alpha(t,\cdot)-\alpha^N(t,\cdot)\|_2+\|\beta(t,\cdot)-\beta^N(t,\cdot)\|_2 +\|\gamma(t,\cdot)-\gamma^N(t,\cdot)\|_2\bigg]=0.$$
This is a straightforward consequence of Theorem \ref{thm-deriv-bd}, for $j=1$. Indeed, since
$$\|\alpha(t,\cdot)-\alpha^N(t,\cdot)\|^2_2=\sum_{k=1}^N\int_{(k-1)/N}^{k/N}|\alpha(t,x)-\alpha(t,(k-1)/N)|^2dx,$$
the bound on $\partial_x\alpha(t,\cdot)$ provided by Theorem \ref{thm-deriv-bd} shows that this quantity tends to vanish as $N\rightarrow\infty$. The same is, obviously, true about the $\beta$ and $\gamma$ terms. Thus, \eqref{main-cv} would follow from
\begin{equation}\label{main-cv-1}
\lim_{N\rightarrow \infty}\bigg[\|a^N(t,\cdot)-\alpha^N(t,\cdot)\|_2+\|b^N(t,\cdot)-\beta^N(t,\cdot)\|_2 +\|c^N(t,\cdot)-\gamma^N(t,\cdot)\|_2\bigg] =0.
\end{equation}
Next, let us prove \eqref{main-cv-1}.
%
Let us define $e^N(t)$ by
\begin{eqnarray*}
&\frac{1}{2}\Big[\left\|a^N(t,\cdot)-\alpha^N(t,\cdot)\right\|^2+\left\|b^N(t,\cdot)-\beta^N(t,\cdot)\right\|^2+\left\|c^N(t,\cdot)-\gamma^N(t,\cdot)\right\|^2\Big]
\end{eqnarray*}
Thus,
\begin{eqnarray*}
e^N(t)&=&\frac{1}{2N}\Big[\left|\vec a^N(t)-\vec\alpha^N(t)\right|^2+\left|\vec b^N(t)-\vec\beta^N(t)\right|^2+\left|\vec c^N(t)-\vec\gamma^N(t)\right|^2\Big]\\
&=&\frac{1}{2N}\left|\vec{u}^N(t)-\vec{v}^N(t)\right|^2.
\end{eqnarray*}
Take the time derivative to see that
\begin{equation}\label{dot-e}
\dot e^N(t)=\frac{1}{N}\Big<\vec{u}^N(t)-\vec{v}^N(t),\frac{d}{dt}\vec{u}^N(t)-\frac{d}{dt}\vec{v}^N(t)\Big>.
\end{equation}
From \eqref{eq1}, \eqref{eqv} and \eqref{dot-e} we obtain
\begin{align}\label{eq:middle}
\dot e^N(t)=\frac{1}{N} &\left<  \vec{u}^N(t)-\vec{v}^N(t),\widetilde{F}^N\left(\vec{u}^N(t)\right)-\widetilde{F}^N\left(\vec{v}^N(t)\right)\right> \\\nonumber
&+ N\left< \vec{u}^N(t)-\vec{v}^N(t),\tilde \Delta^N\left[\vec{u}^N(t)-\vec{v}^N(t)\right]\right>\\\nonumber
&+ \frac{1}{N}\left<\vec{u}^N(t)-\vec{v}^N(t),\vec\varepsilon^N(t)\right>.
\end{align}

For the first term in the right hand side of the above display we use the Mean Value Theorem for vector fields to write
\begin{equation}\label{eqFN}
 \widetilde{F}^N(\vec{u}^N(t))- \widetilde{F}^N(\vec{v}^N(t)) =
 \int_0^1 D\widetilde{F}^N((1-\theta )\vec{u}^N(t)+ \theta \vec{v}^N(t)) \, d\theta\ \vec y(t)
 \end{equation}
where $\vec y(t):=\vec{u}^N(t)- \vec{v}^N(t)$ and  $D\widetilde{F}^N$
denotes the Jacobian matrix of $\widetilde{F}^N$, i.e.
\begin{equation}\label{eq:jac}
D\widetilde{F}^N([\vec a^T, \vec b^T, \vec c^T ]^T) =
\left[
\begin{array}{rrr}
B &A &C \\
B &A &C\\
-B &-A &-C
\end{array}
\right]
\in \mathbb R^{3N\times 3N}
\end{equation}
with generic column vectors $\vec a, \vec b, \vec c\in \mathbb R^N$ and
$A = -k_1 \, \mathrm{diag}(\vec{a}),\ B=-k_1\,\mathrm{diag}(\vec{b})$ and $C=k_{-1}\, I$ ($I$ denotes the identity matrix).
%
%
Equation \eqref{eqFN} yields (we drop the argument $t$ to unburden the notation):
$$\vec y^T[\widetilde{F}^N(\vec{u}^N)-\widetilde{F}^N(\vec{v}^N)] =
\int_0^1 \vec y^T  [D\widetilde{F}^N((1-\theta )\vec{u}^N+\theta \vec{v}^N)] \, \vec y\, d\theta\ .$$
We now fix $[\vec a^T, \vec b^T, \vec c^T]^T:=(1-\theta)\vec{u}^N+\theta \vec{v}^N$ and let $\tilde M$ be a uniform (with respect to $N$, $k$ and $T$) upper bound on the components of $\bf a, \vec b, \vec c$, as per Theorems \ref{the-bound-thm} and \ref{the2}. Then we set the column vector $\vec y=[\vec y^{1T}, \vec y^{2T}, \vec y^{3T}]^T$ (where $\vec y^1, \vec y^2, \vec y^3$ are column vectors in $\mathbb R^N$) and use \eqref{eq:jac} to get
\begin{align}\label{calc_1}
\vec y^T  &[D\widetilde{F}^N((1-\theta )\vec{u}^N+\theta \vec{v}^N)] \, \vec y=[\vec y^{1T}, \vec y^{2T}, \vec y^{3T}]
\left[
\begin{array}{rrr}
B &A &C \\
B &A &C\\
-B &-A &-C
\end{array}
\right]
\left[
\begin{array}{c}
\vec y_1\\
\vec y_2\\
\vec y_3\\
\end{array}
\right]\le \\\nonumber
&\le
[\vec y^{1T}, \vec y^{2T}, \vec y^{3T}]
\left[
\begin{array}{rrr}
k_1\tilde{M}I &k_1\tilde{M}I &k_{-1}I \\
k_1\tilde{M}I &k_1\tilde{M}I &k_{-1}I\\
0&0&0
\end{array}
\right]
\left[
\begin{array}{c}
\vec y^1\\
\vec y^2\\
\vec y^3\\
\end{array}
\right]=\\\nonumber
&=k_1\tilde{M}[(\vec y^{1T}\vec y_1+\vec y^{2T}\vec y_2)+
\vec y^{1T}\vec y^2+\vec y_{2T}\vec y^1)]
+k_{-1}(\vec y^{1T}\vec y^3+\vec y^{2T}\vec y^3)\le\\\nonumber
&\le k1\tilde{M}(|\vec y|^2+\frac{1}{2}|\vec y|^2+\frac{1}{2}|\vec y|^2)+k_{-1}(\frac{1}{2}|\vec y|^2+\frac{1}{2}|\vec y|^2)=\\\nonumber
&=(2k_1\tilde{M}+k_{-1})|\vec y|^2
\end{align}
(this is what is generally known as a bound on the logarithmic norm of the Jacobian).
Thus,
$$\frac{1}{N} \left<  \vec{u}^N(t)-\vec{v}^N(t),\widetilde{F}^N\left(\vec{u}^N(t)\right)-\widetilde{F}^N\left(\vec{v}^N(t)\right)\right> \leq C(\tilde{M})\,e^N(t),$$
where $C(\tilde{M})= 2k_1\tilde{M} + k_1\in \mathbb{R}$ is independent of $N$, $t$ and $T$.\\

The term in the middle of the right hand side of \eqref{eq:middle} is nonpositive because $-\tilde \Delta^N$ is  a nonnegative-definite matrix. Finally, in light of the Cauchy-Schwarz inequality and \eqref{bdd}, the last term in the right hand side of \eqref{eq:middle} is bounded above by $$e^N(t)+\frac{1}{2N}|\vec\varepsilon^N(t)|^2\leq e^N(t)+\frac{3C^2(\delta,T)}{2N^2}.$$
Thus, we have
$$\dot e^N(t)\leq (C(\tilde{M})+1)\,e^N(t)+\frac{3C^2(\delta,T)}{2N^2}\quad\mbox{ for all }\ \delta\leq t\leq T,$$
which implies
\begin{equation}\label{*}
\dot e^N(t) \leq C_1 e^N(t)+ E(N,\delta),
\end{equation}
for some constant $C_1= (C(\tilde{M})+1)$ which is independent of $t$, $N$ and $E(N,\delta):= \frac{3C^2(\delta,T)}{2N^2}.$ Then
 \begin{equation}\label{**}
 \frac{d}{dt}( \exp(-C_1 t) e^N(t))\leq E(N,\delta) \ \exp(-C_1 t),
\end{equation}
Fix $0< t < T$ for given $T>0$, and let $\delta \in (0,t)$. Integrate \eqref{**} from $\delta$ to $t$ to get
\begin{equation*}
e^N(t)\leq \exp( C_1(t- \delta))\ e^N(\delta)+ \frac{E(N,\delta)}{C_1}[\exp( C_1(t- \delta))-1]
\end{equation*}
 and then let $N$ go to infinity to conclude
$$\limsup_{N\rightarrow \infty}e^N(t)\leq \exp( C_1(t- \delta)) \limsup_{N\rightarrow \infty}e^N(\delta)\quad\mbox{ for all } \ t\in[\delta,T].$$
Finally, let $\delta \rightarrow 0^+$ to obtain

\begin{equation}\label{eq:ineq1}
\limsup_{N\rightarrow \infty}e^N(t)\leq C(T) \liminf_{\delta\rightarrow 0^+} \limsup_{N\rightarrow \infty}e^N(\delta),
\end{equation}
where $ C(T)= \exp( C_1 T).$
In view of (\ref{eq:ineq1}) the proof of the theorem is complete once we show that
\begin{equation}
\liminf_{\delta\rightarrow 0^+} \limsup_{N\rightarrow  \infty}  e^N(\delta)  =0,
\end{equation}
which we do next.

Recall that
\begin{align}\label{eq:eNexp}
e^N(\delta) &= \frac{1}{2N}\sum_{k=1}^N \left| \alpha(\delta,kh) - a^N(\delta,kh) \right|^2+  \frac{1}{2N}\sum_{k=1}^N \left| \beta(\delta,kh) - b^N(\delta,kh) \right|^2 + \\\nonumber
&+ \frac{1}{2N}\sum_{k=1}^N \left| \gamma(\delta,kh) - c^N(\delta,kh) \right|^2.
\end{align}
Since
\[
\alpha(\delta,x) = \int_I H_1 (\delta,x,y )\, a_0(y)\,dy + \int_0^\delta \int_I H_1 (s,x,y) f_1(\delta-s,y)\,dy\,ds,
\]
we have
\[
\alpha(\delta,kh) = \int_I H_1 (\delta,kh,y )\, a_0(y)\,dy + \int_0^\delta \int_I H_1 (s,kh,y) f_1(\delta-s,y)\,dy\,ds.
\]
Similarly, using the discrete Heat Kernel $H_1^N$,
\[
a^N(\delta,kh) = \int_I H^N_1 (\delta,kh,y )\, a^N_0(y)\,dy + \int_0^\delta \int_I H^N_1 (s,kh,y) f_1^N(\delta-s,y)\,dy\,ds,
\]
where $a^N_0(y):=\Xint-_{(k-1)h_N}^{kh_N} a_0(u)\,du$ if  $(k-1)h_N \leq y < kh_N$, and \\$f_1^N(s,y):= -k_1a^N(s,y) b^N(s,y) + k_{-1} c^N(s,y).$

Fix $k\in \{1,\hdots,N\}$ to get the estimate
\begin{align}\label{Estim1}
\left| \alpha(\delta,kh) - a^N(\delta,kh) \right| &\leq \bigg| \int_I H_1 (\delta,kh,y )\, a_0(y)\,dy - \int_I H_1^N (\delta,kh,y )\, a_0^N(y)\,dy  \bigg| \\
\nonumber &+  (k_1 M^2 + K_{-1}M) \int_0^\delta \int_I \left|  H_1 (s,kh,y) \right| \,dy\,ds \\
\nonumber &+  (k_1 \tilde{M}^2 + K_{-1}\tilde{M}) \int_0^\delta \int_I \left|  H_1^N (s,kh,y) \right| \,dy\,ds.
\end{align}

$H_1$ is nonnegative and to integrates to 1 in each spatial variable, so
\[
\int_I \left| H_1(\delta,kh,y )\right|\,dy =1 \quad\mbox{for all } \ k=1,\ldots,N,
\]
and it follows that
\begin{align*}
\int_0^\delta \int_I \left| H_1^N (s,kh,y )\right|\,dy\,ds \leq \int_0^\delta \left[ 1 + 2f(4k_As) \right]\,ds = \delta + \frac{1}{2} \int_0^{4k_A\delta} f(s)\,ds.
\end{align*}
Since $f(t)=\sum_{k=1}^{\infty} e^{-k^2 t}$ is positive and integrable on $(0,\infty)$ (see Appendix, Subsection \ref{special_fun}), we have
\[
\lim_{\delta\rightarrow 0^+} \Big[ \delta + \frac{1}{2} \int_0^{4k_A\delta} f(s)\,ds  \Big] = 0.
\]

We have thus obtained bounds on the last two terms in the right hand side of \eqref{Estim1}, depending only on $\delta,$ and not on $k,\, N$. Moreover, these bounds tend to $0$ as $\delta\rightarrow 0^+$. Now focus on the first term in the right hand side of \eqref{Estim1} (call it $T1$). We have
\[
T1\le\int_I H_1(\delta,kh,y ) \left| a_0(y) - a_0^N(y) \right|\,dy + \int_I  \left| H_1(\delta,kh,y ) -  H_1^N (\delta,kh,y ) \right| a_0^N(y)\,dy .
\]
Equation \eqref{BH} in Appendix \ref{special_fun} yields $H_1(\delta,kh,y)\leq 1+2f(4k_A\delta)=:C(\delta)$ for all $y\in I$. Since $a_0^N$ converges in $L^1(I)$ to $a_0$, we may take $N$ sufficiently large so that $ \left\| a_0 - a_0^N \right\|_{L^1(0,1)} \leq \delta/C(\delta)$ (Proposition \ref{l1} in Appendix \ref{a0Ncv}). Since $\left| a_0^N(y) \right| \leq \tilde{M}$ for all $N,\, y$, we get
\begin{align*}
T1 &\leq \delta \int_I H_1(\delta,kh,y )\,dy + \tilde{M}  \int_I  \left| H_1(\delta,kh,y ) -  H_1^N (\delta,kh,y ) \right| \,dy \\
&= \delta + \tilde{M}  \int_I  \left| H_1(\delta,kh,y ) -  H_1^N (\delta,kh,y ) \right| \,dy
\end{align*}
if $N$ is sufficiently large.

But for all $t>0$, $H_1^N (t,\cdot,\cdot )$ converges uniformly to $H_1(t,\cdot,\cdot )$ (see Appendix \ref{cvHNtoH}). Therefore, we have
\[
\tilde{M}  \int_I  \left| H_1(\delta,kh,y ) -  H_1^N (\delta,kh,y ) \right| \,dy \leq \delta
\]
if $N$ is sufficiently large, and so

\begin{align*}
\left| \alpha(\delta,kh) - a^N(\delta,kh) \right| &\leq 2\delta + \delta \left(k_1 M^2 + k_{-1} M\right) + \\
&+\left(k_1 \tilde{M}^2 + k_{-1} \tilde{M}\right) \Big[ \delta + \frac{1}{2} \int_0^{4k_A\delta} f(s)\,ds \Big].
\end{align*}
It is shown similarly that the exact same bound works for $|\beta(\delta,kh) - b^N(\delta,kh)|$ and $| \gamma(\delta,kh) - c^N(\delta,kh)|$, and therefore for sufficiently large $N$ \eqref{eq:eNexp} yields
\[
e^N(\delta) \leq 3  \Big( 2\delta + \delta \left(k_1 M^2 + k_{-1} M\right) + \left(k_1 \tilde{M}^2 + k_{-1} \tilde{M}\right) \Big[ \delta + \frac{1}{2} \int_0^{4k_A\delta} f(s)\,ds \Big] \Big)^2.
\]
The bound $B(\delta)$ above depends on $\delta$ only, and thus
\[
\liminf_{\delta\rightarrow 0^+} \limsup_{N\rightarrow \infty} e^N(\delta)  \leq  \liminf_{\delta\rightarrow 0^+}B(\delta)=\lim_{\delta\rightarrow 0^+}B(\delta)=0.
\]
This concludes the argument for (\ref{eq:ineq1}), and the proof of Theorem \ref{the0}.


\section{Appendix}
\subsection{Heat Kernels}\label{sec:heatKernels}

\paragraph{A. Dirichlet Heat Kernel.}

Let $I:=(0,1)$. Then, $H_D:(0,\infty)\times I\times I \rightarrow \mathbb{R}$ given by
$$H_D(t,x,y) := 2\sum_{j=1}^{\infty} e^{-j^2 \pi^2 t} \sin(j\pi x) \sin(j\pi y)$$
is the Dirichlet Heat Kernel associated to $I$; that is, for any $u_0 \in L^\infty(I)$, the function $u_0:(0, \infty)\times I \rightarrow \mathbb{R}$ given by $$u(t,x)=\int_I H_D(t,x,y)\, u_0(y)\,dy$$ is the unique solution to
\begin{equation}\label{HD}
\begin{cases}
\partial_t u = \partial_x^2 u & \mbox{in $(0,\infty)\times I$}\\
u(\cdot,0) = u(\cdot,1) = 0 & \mbox{in $(0,\infty)$}\\
u(0,\cdot)=u_0(\cdot) & \mbox{in $I$}.
\end{cases}
\end{equation}

\paragraph{Properties of $H_D$}

\begin{enumerate}
\item[$(1)$] $u(t,x)=\int_I H_D(t-\delta,x,y)\, u(\delta,y)\,dy$  for every $t>\delta \geq 0$.
\item[$(2)$] Maximum Principle: $$\max_{\substack{t \in [0,\infty) \\ x\in I}} u(t,x) \leq \max_I u_0.$$
\item [$(3)$] We have $H_D>0$ in $(0,\infty)\times I\times I$, and
$$\int_IH_D(t,x,z)dz\leq 1,\quad \int_IH_D(t,z,y)dz\leq 1\quad\mbox{ for all }x,\ y\in I.$$
\item [$(4)$] It is known that (see, e.g., \cite{Coulhon-Grigor'yan}), there exists a positive constant $C_D>0$ such that
$$H_D(t,x,y)\leq \frac{C_D}{\sqrt{t}}\exp\bigg\{-\frac{(x-y)^2}{8t}\bigg\} \quad\mbox{ for all }t>0,\ x,\ y\in I.$$
\end{enumerate}

Note that $-u$ is also the solution for the initial data $-u_0$. So, by the maximum principle, we also have $$\min_{\substack{t \in [0,+ \infty) \\ x\in I}} u(t,x) \geq \min_I u_0.$$
So, in general, we have
$$\|u(t,\cdot)\|_\infty \leq \|u(\delta,\cdot)\|_\infty \quad\mbox{ for all }\  t >\delta\geq 0.$$

\paragraph{B. Neumann Heat Kernel.}
Let $H:(0,\infty)\times I\times I \rightarrow \mathbb{R}$ given by
\begin{equation}\label{eqHN}
H(t,x,y) := 1+2\sum_{j=1}^{\infty} e^{-j^2 \pi^2 t} \cos(j\pi x) \cos(j\pi y).
\end{equation}
This is the Neumann Heat Kernel associated to $I$; that is, the function (for any given $u_0\in L^\infty(I)$)
\[
u(t,x)=\int_I H(t,x,y)\, u_0(y)\,dy
\]
is the unique solution to
\begin{equation}\label{HN}
\begin{cases}
\partial_t u = \partial_x^2 u & \mbox{in $(0,\infty)\times I$}\\
\partial_x u(\cdot,0)= \partial_x u(\cdot,1)=0 & \mbox{in $(0,\infty)$}\\
u(0,\cdot)=u_0(\cdot) & \mbox{in $I$}.
\end{cases}
\end{equation}

\paragraph{Properties of $H$}:
\begin{enumerate}
\item[$(1')$] $u(t,x)=\int_I H(t-\delta,x,y)\, u(\delta,y)\,dy \quad\mbox{for all } \ t>\delta \geq 0$.
\item[$(2')$] $H\geq 0$ on its domain and $$\int_I H(t,x,y)\,dy = \int_I H(t,x,y)\,dx =1 \quad\mbox{for all } x,\,y\in I, \ \mbox{and all } t>0.$$
\item[$(3')$] From $(1')$ and $(2')$, we also get $$\|u(t,\cdot)\|_\infty \leq \|u(\delta,\cdot)\|_\infty \quad\mbox{for all } \ t >\delta\geq 0.$$
\item[$(4')$] It is known that (see, e.g., \cite{Choulli-Kayser}), there exists a positive constant $C_N>0$ such that
$$H(t,x,y)\leq \frac{C_N}{\sqrt{t}}\exp\bigg\{-\frac{(x-y)^2}{8t}\bigg\}\quad\mbox{ for all }t>0, \ \mbox{and all } x,\ y\in I.$$
\item [$(5')$] $\partial_x H(t,x,y) + \partial_y H_D(t,x,y) = \partial_y H(t,x,y) + \partial_x H_D(t,x,y) =0 \ \mbox{for all } x,\,y\in I, \,\mbox{and }t>0.$
\end{enumerate}

\begin{proposition}\label{deriv-integr}
There exists a positive real number $C$ such that
$$\int_I|\partial_y\tilde H(t,x,y)|dy\leq Ct^{-3/4}\quad\mbox{ for all }(t,x)\in (0,\infty)\times I,$$
where $\tilde H$ is either $H_D$ or $H$.
\end{proposition}
\proof In the proof of Theorem 1.1 \cite{Grigor'yan}, the author shows that if $H_D$ satisfies (on some bounded and open subset $\Omega$ of a smooth, connected, complete noncompact Riemannian manifold $\mathcal{M}$), for every $(t,x)\in (0,T)\times I$, that
\begin{align*}
E_0(t,x) := \int_\Omega|H_D(t,x,y)|^2e^{\frac{(x-y)^2}{4t}}dy \leq \frac{1}{f(t)}
\end{align*}
for some $0<T<\infty$ and some positive $f\in L^1(0,T)$, then
$$E_1(t,x):=\int_\Omega|\partial_yH_D(t,x,y)|^2e^{\frac{(x-y)^2}{4t}}dy\leq \frac{5}{F(t)}\quad\mbox{ for all }(t,x)\in (0,T)\times I,$$
where $F(t):=\int_0^tf(s)ds$. From $(4)$ of the above properties for $H_D$, we see that if we take $\mathcal{M}=\mathbb{R}$ and $\Omega=I$, then we have the desired bound on $E_0$ with $f(t)=t/C_D^2$. We deduce
$$\int_I|\partial_yH_D(t,x,y)|^2e^{\frac{(x-y)^2}{4t}}dy\leq \frac{10\,C_D^2}{t^2} \quad\mbox{ for all }(t,x)\in (0,T)\times I,$$
which, by Cauchy-Schwarz, yields
$$\bigg(\int_I|\partial_yH_D(t,x,y)|dy\bigg)^2\leq \frac{10\,C_D^2}{t^2}\int_Ie^{-\frac{(x-y)^2}{4t}}dy\leq 20\sqrt{\pi}\,C_D^2\,t^{-3/2}.$$
So, in the case $H=H_D$, the statement is proved for $C:=\sqrt{20\sqrt{\pi}}\,C_D$. A careful inspection of the proof of Theorem 1.1 \cite{Grigor'yan} reveals that the same argument works for the Neumann Heat Kernel, so, in light of the property $(4')$ above, we get the desired bound in this case as well.
\endproof

Recall that the solution $(\alpha,\beta,\gamma)$ of \eqref{E1} with Neumann BC satisfies
\begin{align}\label{eqa}
\alpha(t,x)&=\int_I H(t,x,y)\, a_0(y)\,dy + \\\nonumber
&+\int_0^t \int_I  H(t-s,x,y)\,\,[k_{-1} \gamma(s,y)- k_1 \alpha(s,y) \beta(s,y)]\,dy \,ds,
\end{align}
\begin{align}\label{eqb}
\beta(t,x)&=\int_I H(t,x,y)\, b_0(y)\,dy + \\\nonumber
&+\int_0^t \int_I  H(t-s,x,y)\,\,[k_{-1} \gamma(s,y)- k_1 \alpha(s,y) \beta(s,y)]\,dy \,ds,
\end{align}
\begin{align}\label{eqc}
\gamma(t,x)&=\int_I H(t,x,y)\, c_0(y)\,dy +\\\nonumber
&+\int_0^t \int_I  H(t-s,x,y)\,\,[- k_{-1} \gamma(s,y)+ k_1 \alpha(s,y) \beta(s,y)]\,dy \,ds.
\end{align}

  Let us now mimic this representation formula in the discrete case below.\\
\subsection{The Neumann Heat Kernel associated to the discrete case}
Now, we solve the system
\begin{equation}\label{S}
\dot{u}_k^N(t) = \frac{1}{h^2_N}\left( u_{k-1}^N(t) - 2 u_k^N(t) + u_{k+1}^N(t)\right),
\end{equation}
where $k=1,\hdots, N$ and $u_0^N(t):=u_1^N(t),\, u_N^N(t)=:u_{N+1}^N(t)$ (the Neumann BC). We define  $h_N =\frac{1}{N}$.

We can set
\[
U^N(t):=\left[u_1^N(t),\hdots,u_N^N(t)\right]
\]
and rewrite the system \eqref{S} as
\[
\dot{U}^N(t) = \Delta^N U^N(t).
\]
 Note that the matrix $\Delta^N$ has eigenvalues
\[
\lambda_j^N = -4N^2 \sin^2{\frac{(j-1)\pi}{2N}}
\]
and eigenvectors $\vec{\bf v}_j=\left[\nu_{ij}\right]_{1\leq i\leq N}^T$, where
\[
\quad \nu_{ij}=
\displaystyle\begin{cases}
N^{-\frac{1}{2}} & \quad\mbox{if \quad $j=1,\,1\leq i\leq N$}\\
\left(\displaystyle\frac{2}{N}\right)^{\frac{1}{2}} \cos{\displaystyle\frac{(j-1)(i-\frac{1}{2})\pi}{N}} & \quad\mbox{else}.
\end{cases}
\]
If $V^N$ is the $N\times N$ matrix whose columns are $\vec{\bf v}_j,\, j=1,\hdots,N$, then we have
\[
U^N(t)=V^N \exp{\left(D^N t\right)}\left(V^N\right)^T U_0,
\]
where $\exp{\left(D^N t\right)}$ is the diagonal matrix whose diagonal entries are $e^{\lambda_j^N t},\, j=1,\hdots,N$.

Denote by $A^N_{ij}$ the $(ij)^{\mathrm{th}}$ entry in the product $A^N:=V^N\exp{\left(D^N t\right)}\left(V^N\right)^T$.
Define the function $H^N: [0,T]\times I \times I \rightarrow \mathbb{R}$ by
\begin{equation}\nonumber
H^N(t,x,y) = \frac{A^N_{ij}}{h_N} \quad\quad\mbox{if }\quad (i-1)h_N \leq x < i h_N, \quad (j-1)h_N \leq y < j h_N
\end{equation}
for $i,\,j=1,\hdots,N$. Therefore, the solution $U^N(t)$ written as a function $u^N(t,x)$ (defined as $U^N_k(t)$ for $(k-1)h_N\leq x < k h_N$) is given by
\[
u^N(t,x) = \int_I H^N(t,x,y)\, u_0^N(y)\,dy \quad\mbox{for all $(t,x) \in [0,T] \times [0, 1]$.}
\]
From the explicit formulae for the $\vec{\mathbf{v}}_j$'s we compute
\begin{equation}\label{H^N}
H^N(t,x,y) = 1+2 \sum_{k=1}^{N-1} e^{-4N^2\sin^2{\frac{2k\pi}{2N}}}\cos{\frac{k(i - \frac{1}{2})\pi}{N}} \cos{\frac{k(j- \frac{1}{2})\pi}{N}}.
\end{equation}

Going back to the discrete system \eqref{E2}, we have
\begin{equation}\label{eqa^N}
a^N(t,x) = \int_I H^N (t,x,y )\, a^N_0(y)\,dy + \int_0^t \int_I H^N (s,x,y) f_1^N(t-s,y)\,dy\,ds,
\end{equation}
\begin{equation}\label{eqb^N}
b^N(t,x) = \int_I H^N (t,x,y )\, b^N_0(y)\,dy + \int_0^t \int_I H^N (s,x,y) f_2^N(t-s,y)\,dy\,ds,
\end{equation}
\begin{equation}\label{eqc^N}
c^N(t,x) = \int_I H^N (t,x,y )\, c^N_0(y)\,dy + \int_0^t \int_I H^N (s,x,y) f_3^N(t-s,y)\,dy\,ds,
\end{equation}
where $a^N_0(y):=\Xint-_{(k-1)h_N}^{kh_N} a_0(u)\,du $
 $\quad b^N_0(y):=\Xint-_{(k-1)h_N}^{kh_N} b_0(u)\,du$  and $\quad c^N_0(y):=\Xint-_{(k-1)h_N}^{kh_N} c_0(u)\,du$
if  $(k-1)h_N \leq y < kh_N$, and $f_1^N(t,y):= -k_1a^N(t,y) b^N(t,y) + k_{-1} c^N(t,y)$,\\  $f_2^N(t,y):= -k_1a^N(t,y) b^N(t,y) + k_{-1} c^N(t,y),$ and $f_3^N(t,y):= k_1a^N(t,y) b^N(t,y) - k_{-1} c^N(t,y).$

\subsection{Convergence of $H^N(t,x,y)$ to $H(t,x,y)$}\label{cvHNtoH}

Fix $t>0$. Recall that $$H(t,x,y) = 1+ 2\sum_{j=1}^\infty e^{-j^2 \pi^2 t} \cos(j\pi x) \cos(j\pi y) \quad\mbox{for $t\geq 0,\, x,y \in I$},$$
and
\[
H^N(t,x,y) =
\begin{cases}
1+ 2\sum_{j=1}^{N-1} e^{-j^2 \pi^2 t \left(\frac{\sin{\frac{j\pi}{2N}}}{\frac{j\pi}{2N}}\right)^2} \cos{\frac{j(k-\frac{1}{2})\pi}{N}}  \cos{\frac{j(i-\frac{1}{2})\pi}{N}} \\ \mbox{if $t>0,\, \frac{k-1}{N}\leq x <\frac{k}{N},\, \frac{i-1}{N}\leq y <\frac{i}{N},\, i,k=1,\hdots,N$.}
\end{cases}
\]
Of course, in the expression for $H^N$ above, both $k$ and $i$ depend on $N$ and $x,\, y$ (respectively), i.e. $k=k(N,x),\, i=i(N,y)$.

Take an arbitrary $\varepsilon>0$ and fix an integer $m\geq 1$ such that, as the tail of a convergent positive term series, we have
$$ \sum_{j=m+1}^\infty \left(e^{-j^2\pi^2t} + e^{-4 j^2 t}\right) \leq \frac{\varepsilon}{2}.$$
We only consider $N > m$ from now on
and look at
\begin{align*}
\left| H(t,x,y)-H^N(t,x,y)\right| &\leq  2\Big|\sum_{j=1}^m e^{-j^2 \pi^2 t} \cos(j\pi x) \cos(j\pi y) - \\
&-\sum_{j=1}^m e^{-j^2 \pi^2 t \left(\frac{\sin{\frac{j\pi}{2N}}}{\frac{j\pi}{2N}}\right)^2} \cos{\frac{j(k-\frac{1}{2})\pi}{N}}  \cos{\frac{j(i-\frac{1}{2})\pi}{N}}\Big| + \\
&+\sum_{j=m+1}^\infty e^{-j^2\pi^2t} + \sum_{j=m+1}^\infty e^{-4 j^2 t},
\end{align*}
where we have used $\frac{2}{\pi}< \frac{\sin{\frac{j\pi}{2N}}}{\frac{j\pi}{2N}} <1$ for all $1\leq j\leq N$.
With $m$ thus fixed, it remains to show that
\[
\Big|\sum_{j=1}^m \Big[ e^{-j^2 \pi^2 t} \cos(j\pi x) \cos(j\pi y) -  e^{-j^2 \pi^2 t \left(\frac{\sin{\frac{j\pi}{2N}}}{\frac{j\pi}{2N}}\right)^2} \cos{\frac{j(k-\frac{1}{2})\pi}{N}}  \cos{\frac{j(i-\frac{1}{2})\pi}{N}}\Big] \Big| \leq \frac{\varepsilon}{4}
\]
for $N$ sufficiently large. Since $m$ is a fixed positive integer, it is sufficient to prove that for sufficiently large $N$ we have:

\begin{equation}\label{Good}
\Big| e^{-j^2 \pi^2 t} \cos(j\pi x) \cos(j\pi y) -  e^{-j^2 \pi^2 t \left(\frac{\sin{\frac{j\pi}{2N}}}{\frac{j\pi}{2N}}\right)^2} \cos{\frac{j(k-\frac{1}{2})\pi}{N}}  \cos{\frac{j(i-\frac{1}{2})\pi}{N}} \Big| \leq \frac{\varepsilon}{4m}
\end{equation}
for all $j=1,\hdots, m$. So, fix $j\in\{1,\hdots, m\}$. We have
\begin{equation}\label{Eq1}
\lim_{N\rightarrow\infty} \frac{j\pi}{2N}=0\,,\ \mbox{so}\quad \lim_{N\rightarrow\infty} e^{-j^2 \pi^2 t \left(\frac{\sin{\frac{j\pi}{2N}}}{\frac{j\pi}{2N}}\right)^2} = e^{-j^2 \pi^2 t} \quad\mbox{since }\quad \frac{\sin{\frac{j\pi}{2N}}}{\frac{j\pi}{2N}} \xrightarrow[N\rightarrow\infty]{} 1.
\end{equation}
Recall that we also have
$$\frac{k(N,x)-1}{N}\leq x< \frac{k(N,x)}{N}\mbox{ and }\frac{i(N,y)-1}{N}\leq y< \frac{k(N,y)}{N}$$ (where we re-introduced the dependence of $k,\,i$ on $N,\,x,\,y$ to make the point that they vary with $N$ for $x,\,y$ fixed).
It follows that
\[
-\frac{j\pi}{2N}\leq j\pi x -j\pi \frac{k(N,x)-\frac{1}{2}}{N} <\frac{j\pi}{2N} \quad\mbox{and }\quad  -\frac{j\pi}{2N}\leq j\pi y -j\pi \frac{i(N,y)-\frac{1}{2}}{N} <\frac{j\pi}{2N},
\]
i.e., both
\begin{equation}\label{ex2}
\Big| j\pi x -\frac{j (k(N,x)-\frac{1}{2} )\pi}{N} \Big| \leq \frac{j\pi}{2N} \quad\mbox{and }\quad \Big| j\pi y -\frac{j (i(N,y)-\frac{1}{2} )\pi}{N} \Big| \leq \frac{j\pi}{2N}.
\end{equation}
Thus,
\begin{equation}\label{L}\nonumber \tag{$L$}
\lim_{N\rightarrow\infty}\cos{\frac{j(k(N,x)-\frac{1}{2})\pi}{N}} = \cos(j\pi x) \quad\mbox{and }\quad \lim_{N\rightarrow\infty}\cos{\frac{j (i(N,y)-\frac{1}{2} )\pi}{N}} = \cos(j\pi y).
\end{equation}
By \eqref{eq1}, \eqref{L} we get that for each $j\in\{1,\hdots, m\}$, there exists $N(j)$ positive integer such that \eqref{Good} holds for all $N\geq N(j)$. Take $N\geq \max_{j=1,\hdots, m} N(j)$ to conclude
\begin{equation}\label{HNtoH}
\big|H^N(t,x,y)-H(t,x,y)\big|\leq \varepsilon\mbox{ for all }x,\ y\in I.
\end{equation}
Note that $N(j)$ can be chosen independently of $x$ and/or $y$, because the \eqref{L} limits above are approached uniformly with respect to $x,\, y$ (because of \eqref{ex2} and the fact that the cosine function is Lipschitz).

\subsection{A special function}\label{special_fun}

Let $f_n(t):= \sum_{k=1}^n e^{-k^2 t}$ defined on $(0, \infty)$. Clearly, $\left\{f_n\right\}_n$ is an increasing sequence of positive decreasing functions on $(0, \infty)$.

Note that, for every $n\in \mathbb{N}$ and every $T>0$, we have
\begin{align*}
\int_0^T f_n(t)\,dt = {\sum_{k=1}^n \frac{e^{-k^2 t}}{-k^2}\Big|_0^T = \sum_{k=1}^n \frac{1}{k^2}(1-e^{-k^2 T})} \leq {\sum_{k=1}^n \frac{1}{k^2} < \sum_{k=1}^{\infty} \frac{1}{k^2} = \frac{\pi^2}{6}.}
\end{align*}
Thus, $f_n$ is integrable on $(0,\infty)$ for all $n\geq 1$ and $\int_0^{\infty} f_n(t)\,dt \leq \frac{\pi^2}{6}$. By the Monotone Convergence Theorem, the limiting function $f(t):= \sum_{k=1}^{\infty} e^{-k^2 t}$ (which is positive and decreasing) is also integrable on $(0,\infty)$ and $$\int_0^{\infty} f(t)\,dt \leq \frac{\pi^2}{6}.$$

We shall next bound $H$ and $H^N$ in terms of this special function $f$. From \eqref{H^N} we deduce
\begin{align}\label{|H^N|}
\left|H^N(t,x,y)\right| &\leq 1+2 \sum_{k=1}^{N-1} e^{-4N^2\sin^2{\frac{2k\pi}{2N}}} \\ \nonumber &= 1+2 \sum_{k=1}^{N-1} e^{-k^2\pi^2 t \left(\frac{\sin{\frac{2k\pi}{2N}}}{\frac{k\pi}{2N}}\right)^2}.
\end{align}

It is easy to see that $g(x):=\displaystyle\frac{\sin x}{x}$ is positive and decreasing on $[0,\frac{\pi}{2}]$ (at $x=0$ we define $g(0)=1$, obviously). Therefore, we have
\[
\frac{2}{\pi} = \frac{1}{\frac{\pi}{2}}< \frac{\sin{\frac{k\pi}{2N}}}{\frac{k\pi}{2N}}<1 \quad\mbox{for } k=1,\hdots,N-1,
\]
and so
\[
\frac{4}{\pi^2}<\left(\frac{\sin{\frac{k\pi}{2N}}}{\frac{k\pi}{2N}}\right)^2 <1 \quad\mbox{for } k=1,\hdots,N-1.
\]

By \eqref{|H^N|} above, we infer
\begin{align}\label{BH^N}
\left|H^N(t,x,y)\right| &\leq 1+2 \sum_{k=1}^{N-1} e^{-4k^2 t} \\ \nonumber &< 1+2 f(4t).
\end{align}
But in Subsection \ref{cvHNtoH} we proved that for any $t>0$, $H^N(t,\cdot,\cdot)$ converges uniformly to $H(t,\cdot,\cdot)$, so we also get
\begin{align}\label{BH}
\left|H(t,x,y)\right| &\leq 1+2 f(4t).
\end{align}
In particular, for any time interval $[t_1,t_2]$ with $t_2>t_1\geq 0$, we have
\begin{align*}
\int_{t_1}^{t_2}\int_I \left|\tilde H(t,x,y)\right| \,dy\,dt &\leq (t_2-t_1)+2\int_{t_1}^{t_2} f(4 t)\,dt \\
&\leq (t_2-t_1)+\frac{1}{2} \int_0^{\infty} f(t)\,dt \quad\mbox{(change $\pi^2 t \longleftrightarrow t$)}\\
&\leq  (t_2-t_1)+\frac{1}{2} \frac{\pi^2}{6} = (t_2-t_1)+\frac{\pi^2}{12},
\end{align*}
where $\tilde H$ is either $H$ or $H^N$ (for any integer $N\geq2$).

\subsection{Convergence of $a_0^N(x)$ to $a_0$ in $L^p(0,1)$}\label{a0Ncv}

 Let $1\leq p<\infty$ and $a_0 \in L^p(0,1).$ \\
 \begin{proposition}\label{l1} We have
$$ \| a_0^N- a_0 \|_{L^p(0,1)} {\underset{ N\rightarrow \infty} \longrightarrow 0},$$
 where $a^N_0 (x)=  \dashint^{kh}_{(k-1)h} a_0(y)\,dy $ if $(k-1)h\leq x<kh$ for $ k=1,\cdots,N$ and $h = \frac{1}{N}.$
 \end{proposition}
 \proof Take $\delta >0$. Since $C^\infty_c (0,1)$ is dense in $L^p(0,1)$, there exists $\varphi \in C^\infty_c(0,1)$ such that $$\|\varphi-a_0 \|_{L^p(0,1)} \leq \frac{\delta}{3}.$$
 Define $$\varphi ^N(x)= \dashint^{kh}_{(k-1)h} \varphi(y)\,\, dy\mbox{ if } (k-1)h\leq x<kh \mbox{ for } k=1,\cdots, N.$$
 So for sufficient large $N$ we have
 \[
 \|\varphi- \varphi ^N\|^p_{L^p(0,1)} = \sum_{k=1}^N \int^{kh}_{(k-1)h} \bigg| \varphi(x)- \dashint^{kh}_{(k-1)h} \varphi(z)\,\,dz \bigg|^p \,\,dx.
  \]
 Set $\varphi(z_N)= \dashint^{kh}_{(k-1)h} \varphi(z)\,\,dz,$  where $ \frac{k-1}{N} < z_N < \frac{k}{N} $, then by the mean value theorem we have
 $$ \Big| \varphi(x)- \varphi(z_N)\Big| \leq \| \varphi' \|_{\infty} | x-z_N |.$$

 Thus
\begin{align}\label{1}
 \|\varphi-\varphi^N\|^p_{L^p(0,1)} &= \| \varphi' \|^p_{\infty} \sum_{k=1}^N \int^{kh}_{(k-1)h} | x-z_N |^p \,\, dx \leq \\\nonumber
 &\le \| \varphi' \|^p_{\infty} \sum_{k=1}^N \int^{kh}_{(k-1)h} \frac{1}{N^p}\,\, dx = \frac{1}{N^p} \| \varphi' \|^p_{\infty}{\underset{ N\rightarrow \infty} \longrightarrow 0}.
 \end{align}

  Thus $ \|\varphi- \varphi^N\|_{L^p(0,1)}\leq \frac{\delta}{3}$ for sufficient large $N$. \\

 Furthermore,
\begin{align*}
 \| \varphi^N - a_0^N \|^p_{L^p(0,1)}&= \int_I \Big|\varphi^N(x) - a_0^N(x) \Big|^p \,dx \leq \\
 &\le \sum_{k=1}^N \int^{kh}_{(k-1)h} \dashint^{kh}_{(k-1)h} \Big|\varphi(y)- a_0(y)\Big|^p \,dy \,dx \\
 &= \sum_{k=1}^N  \frac{1}{N}\dashint^{kh}_{(k-1)h} \Big|\varphi(y)- a_0(y)\Big|^p \,dy\\
  &= \sum_{k=1}^N  \frac{1}{N} \,  N \int^{kh}_{(k-1)h} \Big|\varphi(y)- a_0(y)\Big|^p \,dy\\
  &= \int_I \Big|\varphi(y)- a_0(y)\Big|^p \,dy,
\end{align*}
Thus,
\begin{equation}\nonumber
 \| \varphi^N - a_0^N \|_{L^p(0,1)} \leq  \| \varphi - a_0 \|_{L^p(0,1)} \leq \frac{\delta}{3}.
\end{equation}
The triangle inequality now yields
\begin{equation}\nonumber
 \|a_0- a_0^N\|_{L^p(0,1)} \leq \|a_0 - \varphi\|_{L^p(0,1)} +  \|\varphi- \varphi ^N\|_{L^p(0,1)} + \|\varphi^N - a_0^N \|_{L^p(0,1)}\leq\delta.
 \end{equation}

\subsection{Multicell networks, complex balanced systems and asymptotic behavior}\label{sec:complexBal}

One of the motivations for this work was the study of asymptotic behavior of complex-balanced reaction-diffusion systems. In the spatially homogeneous case, complex-balanced networks are known to be well behaved, and their study has been central in the field of chemical reaction networks. We briefly, and rather informally, introduce some terminology.

\subsubsection{Chemical reaction networks}

Consider a set of $n$ chemical species with vector of concentrations $x=(x_1, \ldots, x_n)$, and a chemical reaction network (CRN) involving $r$ reactions between these species. Reactions can be viewed formally as arrows between two {\em complexes}, which are formal linear combinations of the species; for example, the reactions (\ref{eq:reac}) considered in this paper are $A+B\to C$ and $C\to A+B$, with complexes $A+B$ and $C$.

Let the system have {\em stoichiometric matrix} $\Gamma$ with rank $r$. Here $\Gamma$ is an $n \times r$ real matrix, and $\Gamma_{ij}$ is the net change in concentration of species $i$ when reaction $j$ occurs. The $j$th column of $\Gamma$ is the {\em reaction vector} for the $j$th reaction. In spatially homogeneous, deterministic, continuous time models, the evolution of the species concentrations is often modeled by mass-action ODEs:
\begin{equation}
\label{genCRN}
\dot x = \Gamma v(x)\,,
\end{equation}
where $v$ is the vector of reaction rates; the rate of each reaction is proportional to the concentrations of reactants. For example, the rate of reaction $A+B\to C$ is $k_1ab$, where $a$ and $b$ denote the concentrations of $A$ and $B$.  The cosets of $\mathrm{im}\,\Gamma$ intersect the nonnegative orthant along {\em stoichiometry classes}. It is easy to see that stoichiometry classes are invariant for (\ref{genCRN}).

A CRN is called {\em complex balanced} \cite{Horn.1972aa} if it admits a positive equilibrium where the net flux at each complex is zero. To make it precise, let $I_{C}, O_{C}\subseteq \{1,\ldots, r\}$ denote the indices of reactions ending and starting at complex $C$. Then $x^*\in{\mathbb R^n_{\ge 0}}$ is a complex balanced equilibrium if for each complex $C$
$$\sum_{j\in I_C} v_j(x^*)\Gamma_j=\sum_{j\in O_C} v_j(x^*)\Gamma_j,$$
where $\Gamma_j$ is the reaction vector of reaction $j$. A CRN is called complex balanced if it admits a positive complex balanced equilibrium, in which case it turns out that all its positive equilibria are complex balanced.  The network $A+B\rightleftharpoons C$ in this paper is trivially complex balanced for any choice of rate constants $k_1$ and $k_2$. More generally, {\em weakly reversible, deficiency zero} networks are complex-balanced for any choice of rate constants \cite{Feinberg.1972aa}. These are networks whose connected components are strongly connected, and for which the number complexes is greater than the number of connected components by rank$\ \Gamma$.

A lot is known about space homogeneous complex balanced systems: they have a unique positive equilibrium in each stoichiometric class, and it is locally asymptotically stable \cite{Horn.1972aa}. A long-standing conjecture states that positive equilibria for complex-balanced systems are in fact {\em globally} asymptotically stable. The reader is referred to \cite{Craciun.2009aa, Anderson.2011sd, Craciun.2013ab, Pantea.2012ss, Gopalkrishnan.2014bb} for partial results towards this conjecture, and  to \cite{Craciun.gac} for a recently announced proof of the general case.

\subsubsection{Multicell reaction networks}
Let $\mathcal R$ be a reaction network with species $X_1,\ldots, X_n$, and let $\Gamma\in{\mathbb R}^{n\times m}$ be its stoichiometric matrix. Fix a positive integer $N$, and let $\bf 1\in{\mathbb R}^N$ denote the column vector of ones.  We let  ${\bf 1}\otimes {\mathcal R}$ define the {\em linear graph  multicell reaction network} \cite{Shapiro.1979aa}, consisting of a collection of $N$ copies ${\mathcal R}^k$ of $\mathcal R$ with species $X_i^k$, $i=\overline{1,n}$, $k=\overline{1,N}$, and additional transport reactions
$$X_i^k\xleftrightharpoons[k_{X_i}]{k_{X_i}} X_i^{k+1},\ i=\overline{1,n},\ k=\overline{1, N-1}.$$

Letting $x_i^k$ denote the concentration of species $X_i^k$, and $x^k=[x_1^k,\ldots, x_n^k]^T$ be the concentration vector of cell $k$,  the ODEs for  ${\bf 1}\otimes {\mathcal R}$ can be written as

\setlength{\tabcolsep}{4pt}
$$\frac{d}{dt}
\left[
\begin{matrix}
x^1\\
\vdots \\
x^N
\end{matrix}
\right]
=
\left[
\begin{tabular}{rrrrrrrrrrrr}
\setlength{\tabcolsep}{4pt}
&$\Gamma$ &-I       &0  &0        &0  &$\ldots$        &0        &0  &0 \\
&0 &I       &$\Gamma$ &-I &0        &0            &$\ldots$        &0  &0\\
&0 &0       &0        &I  &$\Gamma$ &-I &$\ldots$        &0        &0\\
&$\vdots$       &$\vdots$        &$\vdots$  &$\vdots$        &$\vdots$  &$\ddots$ &$\vdots$        &$\vdots$  &$\vdots$\\
&0       &0        &0  &0        &0  &$\ldots$        &$\Gamma$ &-I &0\\
&0       &0        &0  &0        &0  &$\ldots$        &0        &I  &$\Gamma$\\
\end{tabular}
\right]
\left[
\begin{matrix}
v(x^1)\\
w(x^1, x^2)\\
\vdots \\
v(x^{N-1})\\
w(x^{N-1},x^N)\\
v(x^N)\\
\end{matrix}
\right]
$$

The ${nN\times (rN+n(N-1))}$ matrix above is the stoichiometric matrix of ${\bf 1}\otimes {\mathcal R}$, henceforth denoted $\Gamma_{{\bf 1}\otimes {\mathcal R}}$. Here $v(x)$ is the reaction rate vector of $\mathcal R$, and $w(x^k, x^{k+1})=[k_{X_1}(x_1^k-x_2^k), \ldots, k_{X_n}(x_{N-1}^k-x_N^k)]^T$ is the overall transition rate vector between cells $k$ and $k+1$.

Clearly, the conservation laws of ${\bf 1}\otimes {\mathcal R}$ are in direct correspondence with those of $\mathcal R$:
$$\ker\Gamma_{{\bf 1}\otimes {\mathcal R}}^T=\{{\bf 1}\otimes \nu:\ \nu\in\ker\Gamma\},$$
where $\otimes$ denotes the Kronecker product. In other words, for each $\nu\in \ker\Gamma^T,$ $\sum_{i=1}^n \nu_i(x_i^1+x_i^2+\ldots+x_i^N)$ is conserved. Note the similarity with the reaction-diffusion system
$$\partial_t x(t,y)=\Gamma v(x)+ D\Delta x$$
where integrating $\langle \nu, x\rangle$ over the space variable $y$ and using the homogeneous Neumann boundary conditions yields
$\sum_{i=1}^N \int \nu_i x_i(t,y)\ dy = constant.$

Now suppose $y\in\mathbb R_{>0}^n$ is a complex balanced equilibrium of $\mathcal R$. Then it is immediate that ${\bf 1}\otimes y$ is a complex balanced equilibrium of ${\bf 1}\otimes \mathcal R$, and therefore

\begin{proposition}
If $\mathcal R$ is complex balanced, then so is ${\bf 1}\otimes\mathcal R$.
\end{proposition}

This observation has interesting implications: if $\mathcal R$ is complex balanced, then all positive equilibria of ${\bf 1}\otimes\mathcal R$ are asymptotically stable within their compatibility class. This fact, together with the connection made in Theorem \ref{the0} between the reaction-diffusion system (\ref{E1}) and the ODEs corresponding to ${\bf 1}\otimes\mathcal R$, may yield a way of studying the asymptotic behavior of (\ref{E1}), and perhaps of more general classes of complex-balanced systems. This kind of an approach is similar to recent work of Aminzare and Sontag \cite{Aminzare.2014ab, Aminzare.2016aa}, and an alternative to entropy-based techniques
\cite{Desvillettes.2006aa, Desvillettes.2016aa, Fellner.2015aa}. We plan to pursue this line of research in future work.

\vspace{.5cm}

{\bf Acknowledgements.} We thank G. Craciun for encouraging this work, and for informative discussions on multicell reaction networks.

%
%


\begin{thebibliography}{9}

\bibitem{Aminzare.2014ab} Z.~Aminzare and E.~D.~Sontag.
\newblock Synchronization of diffusively-connected nonlinear systems: Results based on contractions with respect to general norms. \newblock {\em IEEE Transactions on Network Science and Engineering}, 1(2):91--106, 2014.

\bibitem{Aminzare.2016aa} Z.~Aminzare and E.~D.~Sontag. \newblock Some remarks on spatial uniformity of solutions of reaction-diffusion {PDE}s. \newblock {\em Nonlinear Analysis: Theory, Methods \& Applications}, 147:125--144, 2016.

\bibitem{Anderson.2011sd} D.~Anderson. \newblock A proof of the global attractor conjecture in the single linkage class case. \newblock {\em SIAM Journal on Applied Mathematics}, 71(4):1487--1508, 2011.

\bibitem{Angeli.2007aa} D.~Angeli, P.~De~Leenheer, and E.~D. Sontag. \newblock A petri net approach to the study of persistence in chemical reaction networks. \newblock {\em Math Biosci}, 210(2):598--618, 2007.

\bibitem{Banaji.2007aa} M.~Banaji, P.~Donnell, and S.~Baigent. \newblock P matrix properties, injectivity, and stability in chemical reaction systems. \newblock {\em SIAM J. Appl. Math}, 67(6):1523--1547, 2007.

\bibitem{Banaji.2016aa} M.~Banaji and C.~Pantea. \newblock Some results on injectivity and multistationarity in chemical reaction networks. \newblock {\em SIAM Journal on Applied Dynamical Systems}, 15(2):807--869, 2016.

\bibitem{WCH} W.~Chen, C.~Li, and E.~Wright. \newblock On A Nonlinear Parabolic System-Modeling Chemical Reactions In Rivers. \newblock {\em Communications On Pure And Applied Analysis}, 4(4):889--899, 2005.

\bibitem{Choulli-Kayser} M.~Choulli and L.~Kayser. \newblock Observations on Gaussian upper bounds for Neumann Heat Kernels, (1991).

\bibitem{ML} M.~Choulli and L.~Kayser. \newblock Observations on Gaussian upper bounds for Neumann Heat Kernels. \newblock {\em Bulletin of the Australian Mathematical Society}, 92(3):429--439, 2015.

\bibitem{Conradi.2007aa} C.~Conradi, D.~Flockerzi, J.~Raisch, and J.~Stelling. \newblock Subnetwork analysis reveals dynamic features of complex (bio) chemical networks.
\newblock {\em PNAS}, 104(49):19175--19180, 2007.

\bibitem{Coulhon-Grigor'yan} T.~Coulhon and A.~Grigor'yan. \newblock Random walks on graphs with regular volume growth. \newblock {\em Geom. Funct. Anal.}, 8(4):656--701, 1998.

\bibitem{Craciun.gac} G.~Craciun. \newblock Toric Differential Inclusions and a Proof of the Global Attractor Conjecture. \newblock {\em arXiv:1501.02860}, 2015.

\bibitem{Craciun.2009aa} G.~Craciun, A.~Dickenstein, A.~Shiu, and B.~Sturmfels. \newblock Toric dynamical systems. \newblock {\em J. Symb. Comp.}, 44(11):1551--1565, 2009.

\bibitem{Craciun.2005aa} G.~Craciun and M.~Feinberg.
\newblock Multiple equilibria in complex chemical reaction networks: I. the injectivity property. \newblock {\em SIAM J. Appl. Math}, 65(5):1526--1546, 2005.

\bibitem{Craciun.2006aa} G.~Craciun and M.~Feinberg.
\newblock Multiple equilibria in complex chemical reaction networks: {II}. the species-reaction graph. \newblock {\em SIAM J. Appl. Math}, 66(4):1321--1338, 2006.

\bibitem{Craciun.2013ab} G.~Craciun, F.~Nazarov, and C.~Pantea. \newblock Persistence and permanence of mass-action and power-law dynamical systems. \newblock {\em SIAM Journal on Applied Mathematics}, 73(1):305--329, 2013.

\bibitem{Desvillettes.2006aa} L.~Desvillettes and K.~Fellner.
\newblock Exponential decay toward equilibrium via entropy methods for reaction-diffusion equations. \newblock {\em Journal of Mathematical Analysis and Applications}, 319(1):157--176, 2006.

\bibitem{Desvillettes.2016aa} L.~Desvillettes, K.~Fellner, and B.~Q.~Tang. \newblock Trend to equilibrium for reaction-diffusion systems arising from complex balanced chemical reaction networks. \newblock {\em arXiv:1604.04536v2}, 2016.

\bibitem{Feinberg.1972aa} M.~Feinberg. \newblock Complex balancing in general kinetic systems. \newblock {\em Archive for Rational Mechanics and Analysis}, 49(3):187--194, 1972.

\bibitem{Fellner.2015aa} K.~Fellner, W.~Prager, and B.~Q.~Tang. \newblock The entropy method for reaction-diffusion systems without detailed balance: first order chemical reaction networks. \newblock {\em arXiv:1504.08221}, 2015.

\bibitem{Gopalkrishnan.2014bb} M.~Gopalkrishnan, E.~Miller, and A.~Shiu. \newblock A geometric approach to the global attractor conjecture. \newblock {\em SIAM J. Appl. Dyn. Sys.}, 13(2):758--797, 2014.

\bibitem{GSW} A.~N.~Gorban, H.~P.~Sargsyan, and H.~A.~Wahab. \newblock Quasichemical models of multicomponent nonlinear diffusion. \newblock {\em Mathematical Modelling of Natural Phenomena}, 6(5):184--262, 2011.

\bibitem{Grigor'yan} A.~Grigor'yan. \newblock Upper bounds of derivatives of the Heat Kernel on an arbitrary complete manifold. \newblock {\em J. Funct. Anal.}, 127(2):363--389, 1995.

\bibitem{Horn.1972ab} F.~Horn. \newblock Necessary and sufficient conditions for complex balancing in chemical kinetics. \newblock {\em Archive for Rational Mechanics and Analysis}, 49(3):172--186, 1972.

\bibitem{Horn.1972aa} F.~Horn and R.~Jackson. \newblock General mass action kinetics. \newblock {\em Archive for Rational Mechanics and Analysis}, 47(2):81--116, 1972.

\bibitem{CE} C. ~Li , E.~S.~Wright. \newblock Modeling chemical reactions in rivers: A three component reaction. \newblock {\em Discrete and Continuous Dynamical Systems}, 7(2):377--384, 2001.

\bibitem{Mincheva.2007aa} M.~Mincheva and M.~Roussel. \newblock Graph-theoretic methods for the analysis of chemical and biochemical networks. {I}. multistability and oscillations in ordinary differential equation models. \newblock {\em Journal of Mathematical Biology}, 55(1):61--86, 2007.

\bibitem{Mincheva.2004aa} M.~Mincheva and D.~Siegel. \newblock Stability of mass action reaction--diffusion systems. \newblock {\em Nonlinear Analysis: Theory, Methods \& Applications}, 56(8):1105--1131, 2004.

\bibitem{Pantea.2012ss} C.~Pantea. \newblock On the persistence and global stability of mass-action systems. \newblock {\em SIAM Journal on Mathematical Analysis}, 44(3):1636--1673, 2012.

\bibitem{Rothe.1984aa} F.~Rothe. \newblock Global solutions of reaction-diffusion systems. \newblock {\em Lecture Notes in Mathematics}. \newblock Vol 1072. \newblock Springer, 1984.

\bibitem{Sadiku} M.~N.~O.~Sadiku and C.~N.~Obiozor. \newblock A simple introduction to the method of lines. \newblock {\em International Journal of Electrical Engineering Education}, 37(3):282--296, 2000.

\bibitem{HJ} A.~H.~Salas, L.~J.~Martinez, and O.~Fernandez. \newblock Reaction-diffusion equations: A chemical application. \newblock {\em Scientia et Technica}, 3(46): 134--137, 2010.

\bibitem{Shapiro.1979aa} A.~Shapiro and F.~Horn. \newblock On the possibility of sustained oscillations, multiple steady states, and asymmetric steady states in multicell reaction systems. \newblock {\em Mathematical Biosciences}, 44(1-2):19--39, 1979.

\bibitem{Siegel.2000aa} D.~Siegel and D.~MacLean.
\newblock Global stability of complex balanced mechanisms.
\newblock {\em J. Math. Chem.}, 27(1):89--110, 2000.

\bibitem{Sontag.2001aa} E.~D.~Sontag. \newblock Structure and stability of certain chemical networks and applications to the kinetic proofreading model of T-cell receptor signal transduction. \newblock {\em IEEE Transactions on Automatic Control}, 46(7):1028--1047, 2001.

\bibitem{SJ} J.~C.~Strikwerda. \newblock Finite Difference Schemes and Partial Differential Equations. \newblock {\em SIAM}, 2004.

\bibitem{Verwer.1984aa} J.~G.~Verwer and J.~M.~Sanz-Serna.
\newblock Convergence of method of lines approximations to partial differential equations. \newblock {\em Computing}, 33(3):297--313, 1984.


\end{thebibliography}
\end{document}